\documentclass[12pt,reqno]{amsart}
\usepackage[T1]{fontenc}
\usepackage[latin9]{inputenc}
\usepackage[letterpaper]{geometry}
\usepackage{color}
\usepackage{amstext}
\usepackage{amsthm}
\usepackage{amssymb}
\usepackage[bookmarks=false,
 breaklinks=false,pdfborder={0 0 1},backref=section,colorlinks=true]
 {hyperref}
\hypersetup{
 breaklinks,linkcolor=blue,citecolor=blue,urlcolor=blue}

\makeatletter

\providecommand{\tabularnewline}{\\}

\numberwithin{equation}{section}
\numberwithin{figure}{section}
\usepackage{enumitem}
\theoremstyle{plain}
\newtheorem{thm}{\protect\theoremname}[section]
\theoremstyle{definition}
\newtheorem{example}[thm]{\protect\examplename}
\theoremstyle{plain}
\newtheorem{lem}[thm]{\protect\lemmaname}
\theoremstyle{plain}
\newtheorem{cor}[thm]{\protect\corollaryname}
\newlist{romaninline}{enumerate*}{1}
\setlist[romaninline]{label=(\roman*)}

\usepackage{graphicx}
\usepackage{color} 
\usepackage{fancyhdr} 
\usepackage{amsfonts} 
\usepackage{latexsym} 
\usepackage{dsfont} 
\usepackage{mathrsfs}
\usepackage{fullpage}
\usepackage{hyperref}
\usepackage{tikz}
\usetikzlibrary{calc, shapes, backgrounds}
\usetikzlibrary{arrows,decorations.pathmorphing,positioning,fit,petri}
\usetikzlibrary{positioning,decorations.pathreplacing,shapes.geometric}
\tikzset{
vtx/.style={inner sep=2pt, outer sep=0pt, circle, fill=black,draw=black}
}

\numberwithin{thm}{section}

\def\COMMENT#1{}
\let\COMMENT=\footnote

\DeclareMathOperator{\mad}{mad}

\def\wt{\widetilde}
\def\mm{~\textrm{mod}^*~}

\makeatother

\providecommand{\corollaryname}{Corollary}
\providecommand{\examplename}{Example}
\providecommand{\lemmaname}{Lemma}
\providecommand{\theoremname}{Theorem}

\begin{document}
\title{Equitable list coloring of sparse graphs}
\author{H. A. Kierstead }
\address{Arizona State University, Tempe, AZ, USA. }
\email{kierstead@asu.edu}
\author{Alexandr Kostochka}
\address{University of Illinois Urbana-{}-Champaign, Urbana, IL 61801 }
\email{kostochk@illinois.edu }
\thanks{Research of AK is supported in part by NSF Grant DMS-2153507 and by NSF RTG Grant
DMS-1937241. }
\author{Zimu Xiang}
\address{University of Illinois Urbana-{}-Champaign, Urbana, IL 61801. }
\email{zimux2@illinois.edu}
\thanks{Research of ZX is supported in part by the Campus Research Board Award RB24000 of
the University of Illinois Urbana-Champaign. }
\begin{abstract}
A proper vertex coloring of a graph is \emph{equitable} if the sizes of all color
classes differ by at most $1$. For a list assignment $L$ of $k$ colors to each
vertex of an $n$-vertex graph $G$, an \emph{equitable} $L$-coloring of $G$ is
a proper coloring of vertices of $G$ from their lists such that no color is used
more than $\lceil n/k\rceil$ times. Call a graph \emph{equitably} $k$-\emph{choosable}
if it has an equitable $L$-coloring for every $k$-list assignment $L$. A graph
$G$ is $(a,b)$-\emph{sparse} if for every $A\subseteq V(G)$, the number of edges
in the subgraph $G[A]$ of $G$ induced by $A$ is at most $a|A|+b$. 

Our first main result is that every $(\frac{7}{6},\frac{1}{3})$-sparse graph with
minimum degree at least $2$ is equitably $3$-colorable and equitably $3$-choosable.
This is sharp. Our second main result is that  every $(\frac{5}{4},\frac{1}{2})$-sparse
graph with minimum degree at least $2$ is equitably $4$-colorable and equitably
$4$-choosable. This is also sharp. 

One of the tools in the proof is the new notion of strongly equitable (SE) list coloring.
This notion is both stronger and more natural than equitable list coloring; and our
upper bounds are for SE list coloring.
\end{abstract}

\maketitle

\section{Introduction}

Let $G=(V,E)$ be a graph, $|G|:=|V|$ and $\|G\|:=|E|$. For $v\in V$, let $N(v)$
be the neighborhood of $v$, $N[v]:=N(v)+v$ and $d(v):=|N(v)|$ be the degree of
$v$. 
 For a subset $A\subseteq V$, let $N(A)=\bigcup_{v\in A} N(v)$ and $N[A]=N(A)\cup A$.
The maximum degree of $G$ is denoted by $\Delta(G)$, and {the minimum degree 
of $G$ 
is denoted by $\delta(G)$.} 
By $V_i(G)$ we denote the set of vertices of degree $i$ in $G$.
The \emph{maximum average degree} of $G$ is defined
by $\mad(G):=\max_{H\subseteq G}$$2\|H\|/|H|$.

\subsection{Equitable coloring.}

In some applications of graph coloring we need not an arbitrary proper coloring, but a coloring in which the color  classes are not too big or are of about the same size. Among such applications, are scheduling in communication systems, construction
timetables, mutual exclusion
scheduling problem, and round-a-clock scheduling, see e.g.~\cite{BEPSW,IL,KW,Mey,SBG,Tu}.
A model for such problems is
 {\em equitable coloring of a graph}---a proper vertex coloring
such that the sizes of every two color classes differ by at most 1. This concept also is useful in studying extremal combinatorial and
probabilistic problems. Alon and F\" uredi~\cite{AF} used it to study existence of some spanning subgraphs in random graphs. Alon and Yuster~\cite{AY} applied results on equitable coloring to the problem of the existence of $H$-factors in dense graphs. Janson and Ruci\' nski~\cite{JR}, Pemmaraju~\cite{Pe} and Janson~\cite{JS} used equitable colorings to
derive deviation bounds for sums of dependent random variables with limited dependency.


The fundamental result on equitable coloring is the Hajnal-Szemer\'{e}di Theorem \cite{HS}
from 1970. It states that every graph $G$ with maximum degree $\Delta(G)=\Delta$
has an equitable $k$-coloring for every $k\geq\Delta+1$.
 The most important open question
on equitable coloring is the Chen-Lih-Wu (CLW) Conjecture~\cite{CLW} from 1994:
Suppose $G$ is a graph with $\Delta(G)=\Delta$; if $G$ has no $(\Delta+1)$-clique,
and when $\Delta$ is odd, no complete bipartite graph with $\Delta$ vertices in
each part, then $G$ has an equitable $\Delta$-coloring. This conjecture is wide
open, but the case $\Delta\leq3$ was proved in \cite{CLW}, 
the case $\Delta=4$ in~\cite{KK12}
and the case $|G|\leq4\Delta$ in~\cite{KK15}.

In a different direction, Bollob\' as and Guy considered equitable coloring of graphs
whose maximum degree is not bounded in terms of $k$. In 1983, they \cite{BG} proved
that a tree $T$ is equitably $3$-colorable if $|T|\geq3\Delta(T)-8$ or $|T|=\Delta(T)-10$
and provided an algorithm for producing the coloring.
This result was extended to all $k\geq 2$ and to all forests by Chen and Lih~\cite{CL} and Miyata,  Tokunaga and Kaneko~\cite{MKT} in 1994.
For a graph $G$ and  $v\in V(G)$, let $\alpha_v=\alpha_v(G)$ denote the size of a maximum independent set in $G$ containing $v$. If a graph $G$ has an equitable $k$-coloring, then by definition, $\alpha_v(G)\geq \lfloor n/k\rfloor$ for each $v\in V(G)$.
Chen and Lih~\cite{CL} and  independently Miyata,  Tokunaga and Kaneko~\cite{MKT} proved (and Chang~\cite{Chang} gave an easier  proof) that this necessary condition is sufficient for forests when $k\geq 3$.

In a similar spirit, Wu and
Wang \cite{WW} proved in 2008 that every planar graph $G$ with minimum degree $\delta(G)\geq2$
and girth $g(G)\geq14$ is equitably $k$-colorable for each $k\geq4$. Moreover,
if the girth of $G$ is at least $26$, then $G$ is equitably $3$-colorable. Luo,
Sereni, Stephens and Yu \cite{LSSY} strengthened this result with the following
theorem.
\begin{thm}[Luo, Sereni, Stephens and Yu]
\label{thm:LuoEtAl}Every planar graph $G$ with $\delta(G)\geq2$ and $g(G)\geq10$
is equitably $k$-colorable for each $k\geq4$, and if the $g(G)\geq14$ then $G$
is equitably $3$-colorable. 
\end{thm}

Luo \emph{et al.} \cite{LSSY} 
mention that their proofs yield the same bounds when they replace
planarity with appropriate bounds on maximal average degree: for graphs $G$ with
$\mad(G)\leq7/3$ when $k=3$ and for graphs with $\mad(G)\leq5/2$ when $k\geq4$,
but they still need the restrictions on the girth and minimum degree.

\subsection{Equitable and strongly equitable list coloring}

A $k$-list assignment for $G$  is a function $L:V\to \binom{\Gamma}{k}$, where $\Gamma$ is a set of \emph{colors}, and $\binom{\Gamma}{k}$ is the set of all $k$-subsets of $\Gamma$. 
Let $\mathcal{L}_k$ be the class of $k$-list assignments for $G$. 
An $L$-coloring of $G$ is a proper coloring $f$ of $G$ with $f(v)\in L(v)$ for all $v\in V$;  $G$ is $L$-colorable if it has an $L$-coloring.  Finally, $G$ is $k$-choosable if it is $L$-colorable for all $L\in \mathcal{L}_k$.

In 2003, Kostochka, Pelsmajer and West (KPW)~\cite{KPW} proposed an equitable version of list coloring. Notice that for any $G$ (with $|G|\ge2k$) there are list assignments $L\in \mathcal{L}_k$, for which there are no $L$-colorings that are equitable colorings in the ordinary sense. For example, if $L(v)\cap L(w)=\emptyset$ and $L(w)=L(x)$ for all $w,x\in V-v$, then any $L$-coloring $G$ has a class of size $1$ and a class of size at least $3$. 
Call a color class $X$ {\em full} if  $X=\lceil|V|/k\rceil$ and \emph{overfull} if $|X|>\lceil|V|/k\rceil$.
For $L\in \mathcal{L}_k$, KPW defined an $L$-coloring to be \emph{list equitable} if it has no overfull classes and graph $G$ to be \emph{equitably $L$-colorable} if it has a list-equitable $L$-coloring. Graph $G$ is \emph{equitably
$k$-choosable} if $G$ has a list-equitable $L$-coloring for all 
$L\in \mathcal{L}_k$.

In 2003, Kostochka\emph{, }Pelsmajer and West (KPW) conjectured that the list
version of the fundamental result holds: every graph $G$ is equitably $k$-choosable
for every $k\geq\Delta(G)+1$. There has been some progress on this conjecture. In 
2013~\cite{KK13} the conjecture was proved for $\Delta\leq7$;
moreover if $|G|\geq\Delta^{3}$ then $G$ is equitably $k$-choosable for every $k\geq\Delta+2$,
and if $|G|\geq4k^{8}$, then $G$ is equitably $(k+1)$-choosable for every $k\ge\max\{\Delta(G),\omega(G)\}$.

It is natural to ask (maybe not conjecture) whether the list version of the CLW Conjecture
is true. In this direction, Dong and Zhang~\cite{DZ18} proved a special case in
2018: 
\begin{thm}[Dong and Zhang]
\label{thm:DongEtAl}If $G$ is a graph with $\mad(G)<3$, then $G$ is equitably
$k$-choosable for $k\ge\max\{\Delta(G),4\}$, and if $\mad(G)<\frac{12}{5}$ then
$G$ is equitably $k$-choosable for $k\ge\max\{\Delta(G),3\}$.
\end{thm}

Recently { we}~\cite{KKX} introduced a stronger and more natural version of list-equitable
 coloring. Consider the $4$-cycle $C$. Let $L$ be the $3$-list assignment for $C$ that assigns all
vertices the set $\{0,1,2\}$, and suppose $f$ is a proper coloring of $C$ using colors  $\{0,1\}$. 
Now $f$ is not an equitable $3$-coloring of $C$ since two classes of $f$ have size $2$, while the third class is empty.
However, $f$ is a list-equitable $L$-coloring of $C$---
no class is overfull.
This is unsatisfactory; the problem is that the definition of list-equitable coloring has allowed too many full classes. An equitable $k$ coloring of a graph on $n$ vertices has exactly $n\mm k$ full classes, where $r=n\mm k$ if  $n=kq+r$, $q,r\in \mathbb{Z}$, and $1\le r\le k$ (so $kq\mm k=k$).  
In~\cite{KKX} we  defined  an $L$-coloring to be \emph{strongly equitable} (SE)
if it has at most $|G|\mm k$ full classes (and no overfull classes). A graph is SE $L$-colorable if it has an SE $L$-coloring. It is
 SE $k$-choosable if it is SE $L$-colorable for
every $L\in\mathcal{L}_{k}$. 

Notice that if we want to extend an SE coloring of
a graph $G_{0}:=G-v$ to an SE $L$-coloring of $G$, then there are at most $d(v)+(|G_{0}|\mod k)$ 
colors that we cannot use for $v$---when $|G_{0}|\mod k=0$, we may have $|G_{0}|\mm k=k$
full classes in $G_{0}$, but at most $1=|G|\mm k$ of these classes will be full
in $G$. In~\cite{KKX}, we proved that every planar graph $G$ is { SE $k$-choosable,
if $k\geq\max\{9,\Delta(G)\}$;} in particular, combining this with our { result in \cite{KK13}} that the KPW Conjecture holds for graphs with $\Delta\le 7$, the KPW Conjecture also holds for every planar graph.

A finer notion than maximum average degree is the notion of $(a,b)$-sparseness.
A multigraph $G$ is $(a,b)$\emph{-sparse}, if for every {$A\subseteq V$ with $|A|\geq 2$},
the number of edges $\|G[A]\|$ induced by $A$ is at most $a|A|+b$. For example,
forests are exactly $(1,-1)$-sparse graphs and planar graphs are $(3,-6)$-sparse.

In this paper, we prove a theorem that unifies, generalizes and strengthens Theorems~\ref{thm:LuoEtAl}
and \ref{thm:DongEtAl} for $k\leq 4$. 
\begin{thm}
\label{thm:main0}Let $G$ be a graph with $\delta(G)\geq2$. If 
\begin{equation}\label{3col}
6\|G[A]\|\leq 7|A|+2\quad\mbox{for each nonempty $A\subseteq V(G)$,} 
\end{equation}
then $G$ SE $3$-choosable, and if 
\begin{equation}\label{4col}
4\|G[A]\|\leq 5|A|+2\quad\mbox{for each nonempty $A\subseteq V(G)$,}
\end{equation}
  then
$G$ is SE $4$-choosable. Moreover, both results are sharp.
\end{thm}

Theorem~\ref{thm:main0} requires a stronger sparseness bound than Theorem~\ref{thm:DongEtAl},
but as in \cite{BG}, allows $\Delta$ to be unbounded and applies to the stronger
parameter, SE choosability. Moreover, it is sharp. Unlike Theorem~\ref{thm:DongEtAl},
it does not assert that $G$ is SE $k$-choosable for larger $k$. We believe that
this is the case, but there are difficulties, and we have put our effort into finding
exact values of sparseness that imply SE $k$-choosability for increasing $k=3,4,\dots$.
We are currently preparing a paper for the harder case $k=5$.

Our results are new for ordinary equitable coloring, as well. They extend the results
mentioned in the remark following Theorem~\ref{thm:LuoEtAl} because they do not
require large girth, and they are tight, even for the finer parameter of $(a,b)$-sparseness. 

The structure of the paper is as follows. In the next section we introduce a number
of notions needed for the proof, state a somewhat more general theorem in a more
convenient language of potentials, setup the proof and indicate its structure. In
Section 3 we show infinitely many sharpness examples for our theorem. In the subsequent
four sections 4--7 we deliver the proof of our more general theorem. In Section
4 we prove results about extending SE list colorings in general graphs. In Section
5 we introduce a minimal counterexample and prove extension results specific to
it. In Section 6 these results are combined to prove lemmas that support the discharging
arguments in Section 7. In Section 8 we conclude the paper with some remarks.

Finally, we should comment on an editorial decision. The two parts of Theorem~\ref{thm:main0}
are proved by very similar arguments, but there are unexpected differences. 
We could prove one part and then address the additional details needed for 
the second, or we could write two stand-alone proofs. The former approach 
invariably shortchanges the second proof, while the latter is repetitious and 
obscures the common thread. We decided to just combine the proofs. Still, 
some readers may prefer at first to concentrate only on one part (we recommend $k=4$), avoiding all references to the other. This is entirely feasible.

\section{Setup}

\subsection{Organization of proof}
It will be easier to prove our result by allowing vertices of degree $1$ in $G$, and modifying inequalities \eqref{3col} and~\eqref{4col}. We also restate them in terms of potentials $\varepsilon_k$ of edges and $\nu_{k}$ of vertices as follows.
Recall that $V_1(G)$ is the set of vertices of degree $1$ in $G$.

Let $\varepsilon_{3}=6$ and
$\varepsilon_{4}=4$.
For each $k\in\{3,4\}$ and subgraph $A\subseteq G$, define  $\nu_{k}:=\varepsilon_{k}+1$,
\begin{equation}\label{nnn}
\rho_{G}^{k}(A)=\varepsilon_{k}\|A\|-\nu_{k}|A|+\frac{\varepsilon_{k}}{2}|A\cap V_{1}(G)|\textrm{ and }\rho_{G}^{k}=\max\{A\subseteq G:\rho_{G}^{k}(A)\}.  
\end{equation}
Note that $A\cap V_{1}(G)$ might not be $V_{1}(G[A])$, since  leaves of $A$ may have neighbors in $V(G)\smallsetminus A$. 

Now  $\rho_{G}^{k}\geq\rho_{G}^{k}(\emptyset)=0$. More importantly,
 $\rho_{G}^{k}$ is
supermodular, i.e., for all $A,B\subseteq V$,
\begin{equation}
\rho_{G}^{k}(A)+\rho_{G}^{k}(B)\leq\rho_{G}^{k}(A\cup B)+\rho_{G}^{k}(A\cap B).\label{eq:supermod}
\end{equation}
Set $\sigma_{G}^{4}:=|V_{1}|\bmod2$ and $\sigma_{G}^{3}=0$. (The particular choice of $\sigma_{G}^{4}$ comes from Example~\ref{exa:4colors}.) The following more
technical result implies our main theorem.
\begin{thm}
\label{thm:Main4} Let $k\in\{3,4\}$. Then all graphs $G$ with $\rho_{G}^{k}\leq2-\sigma_{G}^{k}$
are SE $k$-choosable.
\end{thm}

\subsection{Notation}

\noindent Our notation is motivated by Diestel's   {textbook}~\cite{D}. Let $G=(V,E)$ be
a graph. We assume that $V\cap E=\emptyset$, and notationally treat $G$ as $V\cup E$.
For example, if we know that $v$ is a vertex, then $v\in G$ means the same thing
as $v\in V$. Let $X,Y\subseteq V$ with $X\cap Y=\emptyset$. Denote the set of edges with one end in $X$
and one end in $Y$ by $E(X,Y)$, and put $\|X,Y\|:=|E(X,Y)|$. This notation still
makes sense if $X$ or $Y$ are graphs. 
Set $N(X):=\bigcup_{v\in X}N(v)\smallsetminus X$
and $N[X]:=N(X)\cup X$. Denote a path $P$ with vertices $v_{1},\dots,v_{s}$ and
edges $v_{1}v_{2},\dots,v_{s-1}v_{s}$ by $P=v_{1}\dots v_{s}$. {Set $\mathring{P}:=v_{2}\dots v_{s-1}$,
$\mathring{v_{1}}P:=v_{2}\dots v_{s}$, and $P\mathring{v_{s}}:=v_{1}\dots v_{s-1}$.}
Unlike \cite{D}, we denote the path on  $s$ vertices by $P_s$. 
Denote the cycle $P+v_{s}v_{1}$ by $v_{1}\dots v_{s}v_{1}$.

Define $V_{i}:=V_{i}(G):=\{x\in V:d(x)=i\}$, $V_{i^{+}}:=\bigcup_{j\geq i}V_{j}$,
and $V_{i^{-}}:=\bigcup_{j\leq i}V_{j}$.
Given a graph $G$ and a color $\alpha$ used by some coloring $f$ of $G$, let $\wt\alpha=\wt\alpha(f)$ be the set of vertices of $G$ colored with $\alpha$ by $f$.
Let $T(f)$ be the set of full colors regarding $f$.
For a subset $S\subseteq V(G)$, set $f(S):=\{f(v):v\in S\}$.
\section{Sharpness examples}
{ In this section, we present infinitely many  examples of graphs with $\rho^k=3-\sigma^k$ that admit no equitable $k$-coloring, where $k=3,4$. This immediately implies they are not SE $3$-choosable, and
so  Theorem~\ref{thm:Main4} is sharp. 
Moreover, infinitely many of these examples have order divisible by $k$, so they are also not equitably $k$-choosable.}

\begin{example}
For all $n\in\mathbb{N}$ and $l=0,\dots,5$, there are triangle-free outerplanar graphs $G_{n,l}$
with $|G_{n,l}|>n$, $|V_{1}(G_{n,l})|=l$ and $\rho_{G_{n,l}}^{3}=3$ that admit no equitable
$3$-coloring.
\end{example}

\begin{proof}
For $l\in\{0,\dots,5\}$ and $n\in\mathbb{N}$, graph $G_{n,l}$ has a special vertex $x$ and consists of $n$ copies of $C_7$, $l$ copies of $K_2$ and $5-l$ copies of $C_5$ such that each of these copies contains $x$ and all of them are otherwise disjoint (See Figure~\ref{fig:k=00003D3}.). Note that each $G_{n,l}$ is outerplanar.

Then $|G_{n,l}|=1+6n+l+4(5-l)=21+6n-3l$ and $|V_{1}(G_{n,l})|=l$.

If $G_{n,l}$ would have an equitable $3$-coloring, then every color class would have $7-l+2n$ vertices. But the color class $C_x$ of $x$ apart from $x$ itself has at most $2$ vertices in each copy of $C_7$, at most one vertex in each copy of $C_5$ and no other vertices, so that
$|C_x|\leq 1+2n+(5-l)<7-l+2n$. Thus $G_{n,l}$  has no equitable $3$-coloring. 

Finally, we show that $\rho_{G_{n,l}}^{3}=3$, for all $n\in\mathbb{N}$ and $l=0,\dots,5$.
Consider a subgraph $A\subseteq G_{n,l}$ with $\rho_{G_{n,l}}^{3}(A)$ maximum. If $x\notin V(A)$, then each component of $A$ is either a singleton in $V_1(G)$ or a path disjoint from $V_1(A)$. In both cases, the potential of each component is negative. Thus $x\in V(A)$. The maximality of 
$\rho_{G_{n,l}}^{3}(A)$ then implies that all $l$ vertices in $V_1(G)$ are in $A$, each of the $(5-l)$ copies of $C_5$ is in $A$ and for each copy $B$ of $C_7$, either the whole $B$ is in $A$ or only $x$. This means that the potential of $A$ is exactly $-7+l(6-4)+(5-l)(6(5)-7(4))=3$.
\end{proof}

\begin{figure}
\begin{tikzpicture}[scale =1] 
\def \vt{circle (1.5pt) [fill]} \def \mg{cyan!50!red!20!white} \def \mc{red!20!yellow!20!white} \def \mpg{cyan!20!green!20!white} \foreach \j in {-10,35,80,125,170} { \draw (\j:2cm)--(\j+20:2cm)--(0:0cm)--cycle; 
\draw (\j:2cm) \vt; \draw (\j:1cm) \vt; \draw (\j+20:2cm) \vt; \draw (\j+20:1cm) \vt;} 
\draw (0,0) \vt; 
\draw[red] (0,0)--(.5,-.5)--(.5,-1)--(.25,-1.5)--(-.25,-1.5)--(-.5,-1)--(-.5,-.5)--cycle; 
\draw[red] (.5,-.5) \vt; 
\draw[red] (.5,-1) \vt;
\draw[red] (-.5,-1) \vt; 
\draw[red] (-.5,-.5) \vt; 
\draw[red] (.25,-1.5) \vt;
\draw[red] (-.25,-1.5) \vt;
\draw (0,0) \vt; 
\end{tikzpicture}
\hspace{2cm}
\begin{tikzpicture}[scale =1] 
\def \vt{circle (1.5pt) [fill]} \def \mg{cyan!50!red!20!white} \def \mc{red!20!yellow!20!white} \def \mpg{cyan!20!green!20!white} \foreach \j in {80,125,170} { \draw (\j:2cm)--(\j+20:2cm)--(0:0cm)--cycle; 
\draw (\j:2cm) \vt; \draw (\j:1cm) \vt; \draw (\j+20:2cm) \vt; \draw (\j+20:1cm) \vt;} 
\foreach \j in {0,37.5} { \draw (\j:2cm)--(0:0cm); 
\draw (\j:2cm) \vt; } 
\draw (0,0) \vt; 
\draw[red] (0,0)--(.5,-.5)--(.5,-1)--(.25,-1.5)--(-.25,-1.5)--(-.5,-1)--(-.5,-.5)--cycle; 
\draw[red] (.5,-.5) \vt; 
\draw[red] (.5,-1) \vt;
\draw[red] (-.5,-1) \vt; 
\draw[red] (-.5,-.5) \vt; 
\draw[red] (.25,-1.5) \vt;
\draw[red] (-.25,-1.5) \vt;
\draw (0,0) \vt; 
\end{tikzpicture}
\hspace{2cm}
\begin{tikzpicture}[scale =1] 
\def \vt{circle (1.5pt) [fill]} \def \mg{cyan!50!red!20!white} \def \mc{red!20!yellow!20!white} \def \mpg{cyan!20!green!20!white} \foreach \j in {30,60,90,120,150} { \draw (\j:2cm)--(0:0cm); 
\draw (\j:2cm) \vt; } 
\draw (0,0) \vt; 
\draw[red] (0,0)--(.5,-.5)--(.5,-1)--(.25,-1.5)--(-.25,-1.5)--(-.5,-1)--(-.5,-.5)--cycle; 
\draw[red] (.5,-.5) \vt; 
\draw[red] (.5,-1) \vt;
\draw[red] (-.5,-1) \vt; 
\draw[red] (-.5,-.5) \vt; 
\draw[red] (.25,-1.5) \vt;
\draw[red] (-.25,-1.5) \vt;
\draw (0,0) \vt; 
\end{tikzpicture}

\caption{\protect\label{fig:k=00003D3}$\rho_{G}^{3}=3$ and $G$ cannot be equitably $3$
colored, where $G=G_{1,0}$ on the left, $G=G_{1,2}$ in the middle and $G=G_{1,5}$ on the right.}
\end{figure}

\begin{example}\label{exa:4colors}
For all $n\in\mathbb{N}$ and $l=0,\dots,7$, there are outerplanar graphs $G_{n,l}$
with $|G_{n,l}|>n$, $|V_{1}(G_{n,l})|=l$ and $\rho_{G}^{4}=3-\sigma^4(G_{n,l})$
that admit no equitable $4$-coloring.
\end{example}

\begin{proof} Let $n\in\mathbb{N}$,  $l\in\{0,\dots,7\}$, $l':=l\pmod 2$ and
$r=\frac{8-l-l'}{2}$. Graph $G_{n,l}$ has a special vertex $x$ and consists of $n$ copies of $C_5$, $l$ copies of $K_2$ and $r$ copies of $C_3$ such that each of these copies contains $x$ and all of them are otherwise disjoint (See Figure~\ref{fig:1}.).
Note that each $G_{n,l}$ is outerplanar.
Then $|G_{n,l}|=1+4n+l+2\frac{8-l+l'}{2}=9+4n-l'$ and $|V_{1}(G_{n,l})|=l$.

If $G_{n,l}$ would have an equitable $3$-coloring, then every color class would have at least $2+n$ vertices. But the color class $C_x$ of $x$ apart from $x$ itself has at most $1$ vertex in each copy of $C_5$ and no other vertices, so that
$|C_x|\leq 1+n$. Thus $G_{n,l}$  has no equitable $4$-coloring. 

Finally, we show that $\rho_{G}^{4}\leq 3-\sigma_{G_{n,l}}^{4}$, for all $n\in\mathbb{N}$
and $l=0,\dots,7$.
Consider a subgraph $A\subseteq G_{n,l}$ with $\rho_{G_{n,l}}^{4}(A)$ maximum. If $x\notin V(A)$, then each component of $A$ is either a singleton in $V_1(G)$ or a path disjoint from $V_1(A)$. In both cases, the potential of each component is negative. Thus $x\in V(A)$. The maximality of 
$\rho_{G_{n,l}}^{3}(A)$ then implies that all $l$ vertices in $V_1(G)$ are in $A$, each of the $r$ copies of $C_3$ is in $A$ and for each copy $B$ of $C_5$, either the whole $B$ is in $A$ or only $x$. This means that the potential of $A$ is exactly $-5+l(4-3)+r(4(3)-5(2))=-5+l+2\frac{8-l-l'}{2}=3-l'=3-\sigma^4(G_{n,l})$.
\end{proof}

\begin{figure}
\begin{tikzpicture}[scale =1]
\def \vt{circle (1.5pt) [fill]} \def \mg{cyan!50!red!20!white} \def \mc{red!20!yellow!20!white}  \def \mpg{cyan!20!green!20!white} \foreach \j in {20,60,100,140} { \draw (\j:2cm)--(\j+20:2cm)--(0:0cm)--cycle; 
\draw (\j:2cm) \vt; \draw (\j+20:2cm) \vt; } 
\draw (0,0) \vt; 

\draw[red] (0,0)--(.5,-.5)--(.5,-1)--(-.5,-1)--(-.5,-.5)--cycle; 
\draw[red] (.5,-.5) \vt; 
\draw[red] (.5,-1) \vt;
\draw[red] (-.5,-1) \vt; 
\draw[red] (-.5,-.5) \vt; 
\draw (0,0) \vt;
\end{tikzpicture}
\hspace{2cm}
\begin{tikzpicture}[scale =1]
\def \vt{circle (1.5pt) [fill]} \def \mg{cyan!50!red!20!white} \def \mc{red!20!yellow!20!white}  \def \mpg{cyan!20!green!20!white} 
\foreach \j in {100,140} { \draw (\j:2cm)--(\j+20:2cm)--(0:0cm)--cycle; 
\draw (\j:2cm) \vt; \draw (\j+20:2cm) \vt; } 
\draw[red] (0,0)--(.5,-.5)--(.5,-1)--(-.5,-1)--(-.5,-.5)--cycle; 
\draw[red] (.5,-.5) \vt; 
\draw[red] (.5,-1) \vt;
\draw[red] (-.5,-1) \vt; 
\draw[red] (-.5,-.5) \vt; 
\draw (0,0) \vt;

\def \vt{circle (1.5pt) [fill]} \def \mg{cyan!50!red!20!white} \def \mc{red!20!yellow!20!white}  \def \mpg{cyan!20!green!20!white} \foreach \j in {30,50,70} { \draw (\j:2cm)--(0:0cm); 
\draw (\j:2cm) \vt; } 
\draw (0,0) \vt; 

\draw[red] (0,0)--(.5,-.5)--(.5,-1)--(-.5,-1)--(-.5,-.5)--cycle; 
\draw[red] (.5,-.5) \vt; 
\draw[red] (.5,-1) \vt;
\draw[red] (-.5,-1) \vt; 
\draw[red] (-.5,-.5) \vt; 
\draw (0,0) \vt;
\end{tikzpicture} 
\hspace{2cm}
\begin{tikzpicture}[scale =1]
\def \vt{circle (1.5pt) [fill]} \def \mg{cyan!50!red!20!white} \def \mc{red!20!yellow!20!white}  \def \mpg{cyan!20!green!20!white} 

\draw[red] (0,0)--(.5,-.5)--(.5,-1)--(-.5,-1)--(-.5,-.5)--cycle; 
\draw[red] (.5,-.5) \vt; 
\draw[red] (.5,-1) \vt;
\draw[red] (-.5,-1) \vt; 
\draw[red] (-.5,-.5) \vt; 
\draw (0,0) \vt;

\def \vt{circle (1.5pt) [fill]} \def \mg{cyan!50!red!20!white} \def \mc{red!20!yellow!20!white}  \def \mpg{cyan!20!green!20!white} \foreach \j in {30,50,70,90,110,130,150} { \draw (\j:2cm)--(0:0cm); 
\draw (\j:2cm) \vt; } 
\draw (0,0) \vt; 

\draw[red] (0,0)--(.5,-.5)--(.5,-1)--(-.5,-1)--(-.5,-.5)--cycle; 
\draw[red] (.5,-.5) \vt; 
\draw[red] (.5,-1) \vt;
\draw[red] (-.5,-1) \vt; 
\draw[red] (-.5,-.5) \vt; 
\draw (0,0) \vt;
\end{tikzpicture} 

\caption{\protect\label{fig:1}$\rho_{G}^{4}=3-\sigma_{G}^{4}$ and $G$ cannot be equitably
$4$ colored, where $G=G_{1,0}$ on the left, $G=G_{1,3}$ in the middle and $G=G_{1,7}$ on the right.}
\end{figure}

\section{Safe bugs in general graphs}

In this section, we consider arbitrary graphs $G$, integers $k\geq3$ and list assignments
$L\in\mathcal{L}_{k}.$ Call a subgraph $S\subseteq G$ \emph{safe} in $G$ if every
SE $L$-coloring of $G_{0}:=G-S$ can be extended to an SE $L$-coloring of $G$.
The safety of $S$ does not imply that $G_{0}$ is SE $L$-colorable. For disjoint
induced subgraphs $S_{0}$ and $S_{1}$, 
\begin{equation}
\text{if \ensuremath{S_{0}} is safe in~\ensuremath{G-S_{1}} and \ensuremath{S_{1}} is safe in \ensuremath{G}, then \ensuremath{S_{0}\cup S_{1}} is safe in \ensuremath{G}. }\label{eq:trans}
\end{equation}

The following is a modification of Lemma 1.1 in~\cite{MP} by Pelsmajer.
\begin{lem} 
\label{lem:P+}Suppose $S=G[\{v_{1},\dots,v_{s}\}]$ with $|S|=s\leq k$. Then $S$
is safe in $G$ if
\begin{equation}
\|v_{i},G-S\|<i\text{ for all }i\in[s].\label{eq:safe-1}
\end{equation}
\end{lem}

\begin{proof}{ Suppose $G_{0}:=G-S$ has an SE $L$-coloring $f_{0}$. Let $m:=|G_{0}|\mm k$. By
definition, $|T(f_{0})|\le m$. Order the colors so that $|f^{-1}(\alpha)|<|f^{-1}(\beta)|$ implies $\alpha<\beta$. Extend $f_{0}$ to a coloring $f$ of $G$ by coloring for $i=1,2,\ldots,s$
\, vertex $v_{i}$ with the least color in $L(v_{i})$ that is not used on any vertex in $(N(v_{i})\smallsetminus S)\cup \{v_{i+1},\dots,v_{s}\}$.  
There are at most $(i-1)+(s-i)=s-1$ forbidden colors, 
so $v_{i}$ is colored with one of the $s$-smallest colors, and  (i) the $k-s$ largest colors are not used on $S$. Now $f$ is proper. 
Also, (ii) no color class has gained more than one vertex. 
If $s> k-m$ then $\lceil\frac{n}{k}\rceil=\lceil\frac{n-s}{k}\rceil+1$.
Thus by (ii)  $T(f)\subseteq T(f_{0})$, and by (i) $|T(f)|\le m-(k-s)=n\mm k$, so $f$ is SE. 
  Else $s\le k-m$. By (i) each $v_{i}$  is colored with a color $\alpha\notin T(f_{0})$, so $f$ is equitable. As $|T(f)|\le m+s\le n\mm k$, $f$ is SE. Thus $S$ is safe.}
\end{proof}
A path or cycle $P:=r\dots a\subseteq G$ is a \emph{thread} if $r\in V_{3^{+}}(G)$
and $P\cap V_{3^{+}}\subseteq\{r,a\}$; $P$ is \emph{plain} if {$a\in V_{3^{+}}-r$},
\emph{loose} if $a\in V_{1}(G)$ and \emph{closed} if $a=r$. Depending on the context,
we may refer to $P$ as an $r,a$-thread if $|P\cap V_{3^{+}}|=2$, or as a $t$-thread
where $t=|P-V_{3^{+}}|$. 

A \emph{bug }$B\subseteq G$ is a connected subgraph with a vertex $r$, called the
\emph{root} of $B$, such that $B\cap V_{3^{+}}(G)\subseteq\{r\}$. Figure~\ref{fig:Bug} shows some bugs. A vertex $v\in B-r$ 
is a \emph{leg} vertex if it is in a plain thread with an end $r$; it is a \emph{body}
vertex if it is in a loose or closed thread ending at $r$. 
 Let $B(r)=B_{G}(r)$ be the maximal bug with
root $r$, that is the bug containing every bug containing $r$.
For a subgraph $H\subseteq G$, a vertex
$v\in H$ is \emph{hidden} if $N_{G}(v)\subseteq H$. 

\begin{figure}
\begin{tikzpicture}
\def \vt{circle (1.5pt) [fill]} 
\draw (0,0)--(-1,1);
\draw (0,0)--(1,1);
\draw(-1,1)--(1,1);
\draw(0,0)--(-1,-1);
\draw (0,0)--(1,-1);
\draw(-1,-1)--(-1.25,-1.25);
\draw(1,-1)--(1.25,-1.25);
\draw (0,0) \vt; 
\draw (0.4,0) node(r) {$r$};
\draw (-1,1) \vt;
\draw (-1.4,1) node(v3) {$v_3$};
\draw (1,1) \vt;
\draw (1.4,1) node(v4) {$v_4$};
\draw (-1,-1) \vt;
\draw (-1.4,-1) node(v1) {$v_1$};
\draw (1,-1) \vt;
\draw (1.4,-1) node(v2) {$v_2$};
\draw (0,-2) node(a) {(a)};
\end{tikzpicture}
\hspace{.5cm}
\begin{tikzpicture}
\def \vt{circle (1.5pt) [fill]}
\draw (0,1)--(0,0)--(0,-1);
\draw (0,0)--(1,-1);
\draw (0,0)--(-1,-1);
\draw (1.25,-1.25)--(1,-1);
\draw (-1.25,-1.25)--(-1,-1);
\draw (0,-1)--(0,-1.25);
\draw (0,0) \vt; 
\draw (0.4,0) node(r) {$r$};
\draw (0,1) \vt;
\draw (0.4,1) node(v4) {$v_4$};
\draw (1,-1) \vt;
\draw (1.4,-1) node(v3) {$v_3$};
\draw (-1,-1) \vt;
\draw (-1.4,-1) node(v1) {$v_1$};
\draw (0,-1) \vt;
\draw (-0.4,-1)  node(v2) {$v_2$};
\draw (0,-2) node(b) {(b)};
\end{tikzpicture}
\hspace{.5cm}
\begin{tikzpicture}
\def \vt{circle (1.5pt) [fill]}
\draw (0,1)--(0,-.5);
\draw (0,1)--(-1.5,-.5);
\draw (0,1)--(1.5,-.5);
\draw (0,1) \vt;
\draw (0.4,1) node(r) {$r$};
\draw (0,0) \vt;
\draw (-0.4,0) node(v2) {$v_2$};
\draw (-1,0) \vt;
\draw (-1.4,0) node(v1) {$v_1$};
\draw (0.6,0.4) \vt;
\draw (1,0.4) node(x) {$x$};
\draw (1.2,-.2) \vt;
\draw (1.6,-.2) node(y) {$y$};
\draw (0,-1.5) node(c) {(c)};
\end{tikzpicture}
\hspace{.5cm}
\begin{tikzpicture}
\def \vt{circle (1.5pt) [fill]}
\draw (0,1)--(0.5,-.5);
\draw (0,1)--(-0.45,-.5);
\draw (0,1)--(1.5,-.5);
\draw (0,1)--(-1.5,-.5);
\draw (0,1) \vt;
\draw (0.4,1) node(r) {$r$};
\draw (0.333,0) \vt;
\draw (0.6,0) node(v2) {$v_2$};
\draw (-0.3,0) \vt;
\draw (-.6,0) node(v1) {$v_1$};
\draw (0.6,0.4) \vt;
\draw (1,0.4) node(x) {$x$};
\draw (1.2,-.2) \vt;
\draw (1.6,-.2) node(y) {$y$};
\draw (0,-1.5) node(d) {(d)};
\end{tikzpicture}
\caption{\label{fig:Bug} { The figure shows four maximal bugs. In (a), 
$v_1,v_2$ are leg vertices, $v_3,v_4$ are body vertices, and $v_3,v_4,r$ are hidden; 
in (b), $v_4$ is a body vertex and  $v_4,r$ are hidden; in (c), 
$v_1,v_2,x,y$ are leg vertices, and  $x,r$ are hidden; in (d), $v_1,v_2,x,y$ are leg vertices, and  $x$, but not $r$, is hidden.}} 
\end{figure}

\begin{cor}
\label{cor:safe}{Let $B\subseteq G$ be a bug with root $r$. If $\|r,G-B\|<|B|\leq k$ and
 $B$ has a hidden vertex $y_{1}$, then $B$ is safe in $G$.}
\end{cor}

\begin{proof}
Set $G_{0}:=G-B$.  Suppose $x\in B-r$. Then $\|x,G_{0}\|=d_G(x)-d_{B}(x)$,   $d_{G}(x)\leq2$
and $d_{B}(x)\geq1$, so $\|x,G_{0}\|\leq1$. As $y_1$ is hidden and $\|r,G_0\|<|B|$, the ordering $y_{1},\dots,y_{|B|}$
with $r=y_{|B|}$,  satisfies (\ref{eq:safe-1}). Thus, $B$ is safe
in $G$ by Lemma~\ref{lem:P+}. 
\end{proof}

\begin{cor}
\label{cor:bsafe}Let $B$ be a bug with root $r\in V_{k^{-}}$ and $k+1\leq|B|\leq2k$.
If $B$ has distinct nonadjacent hidden vertices, then $B$ is safe in $G$. 
\end{cor}

\begin{proof}
Suppose there are nonadjacent hidden vertices $x,y\in B$ and a set $B_{2}\subseteq B$
such that
\begin{romaninline}
\item $r=x$ or $r\in N(y)$,
\item $N[y]\subseteq B_{2}\subseteq B-x$,
\item $|B_{2}|=k$, and
\item some $z\in B_{2}\smallsetminus\{r,y\}$ satisfies $\|z,G - B_2\|\leq 1$.
\end{romaninline}
Set $B_{1}:=B-B_{2}$. Now using Lemma~\ref{lem:P+} twice, $B_{1}$ is safe in
$G_{1}:=G-B_{2}$, since $x$ is hidden, $\|v,V\smallsetminus B\|\leq1$ for all
$v\in B-r$ and $r\notin B_{1}-x$ by (i) and (ii); and $B_{2}$ is safe in $G$,
since $y$ is hidden in $B_{2}$ by (ii), some $z\in B_{2}\smallsetminus\{r,y\}$ satisfies $\|z,G_{1}\|\leq1$
by (iv), and, as $r\in V_{k^-}$, if $r\in B_{2}$ then $\|r,G_{1}\|\leq_{\text{(i,ii)}}k-1<_{\text{(iii)}}|B_{2}|$, and all other vertices $v$ in $B_2$ have $\|v,G_1\| \leq 2$ by the definition of $B$.
So $B$ is safe in $G$ by (\ref{eq:trans}).

{It remains to pick $x,y,B_{2}$ as above. First, we choose distinct nonadjacent hidden
$x$ and $y$ satisfying (i). If $r$ is hidden and there is a hidden $v\in B\smallsetminus N[r]$,
then set $x:=r$ and $y:=v$. Else, pick any distinct nonadjacent hidden vertices
$x$ and $y$ with $d(y)\geq d(x)$.  Anyway, $N[y]\subseteq B-x$, and by hypothesis,
$|B|\geq k+1$ and the fact that $k\ge 3$, so there is $B_{2}$ satisfying (ii) and (iii). If $Z:=N(y)-r\ne\emptyset$
then (iv) holds for any $z\in Z$. If $Z=\emptyset$ then $N(y)=\{r\}$ and $d(x)\le d(y)=1$.
Now $B-x$ is connected, so we can choose $B_{2}$ to be connected. Then (iv) holds
for any $z\in B_{2}\smallsetminus\{r,y\}$.}
\end{proof}
\begin{cor}
\label{cor:d(v)=00003D1}If $B$ is a bug with root $r\in V_{1^{-}}$. Then $B$
is safe in $G$ if $|B|\geq2$ or $\|B,G-B\|=0$ or $|G|\ne0\mod k$. 
\end{cor}

\begin{proof}
Set $G_{0}:=G-B$. As $d(r)=1$, $B$ is a path $P=v_{1}\dots v_{t}$ with $r=v_{1}$
and $N(B)\subseteq N(v_{t})$. If $2\leq t\leq k$ or $t=1$ and $\|B,G-B\|=0$, then
$r$ is hidden, so $B$ is safe in $G$ by Corollary~\ref{cor:safe}. Arguing by
induction, suppose $t>k$. By induction, $v_{3}Pv_{t}:=v_3\dots v_t$ is safe in $G_{1}:=G-v_{1}v_{2}$
and $v_{1}v_{2}$ is safe in $G$ by Lemma~\ref{lem:P+}. By (\ref{eq:trans}), $P$ is safe in $G$. 

Otherwise $t=1$ and $r$ has a neighbor $a\notin B$. If $B$ is not safe, then some
SE $L$-coloring $f$ of $G_{0}$ cannot be extended to $G$. Thus $L(r)-f(a)\subseteq T(f)$
and $|G_{0}|\ne0\bmod k$. Thus $|T(f)|=k-1$, $|G_{0}|=k-1\bmod k$, and $|G|=0\bmod k$.
\end{proof}
\begin{cor}
\label{cor:d(v)=00003D2} If $B$ is a bug with root $r\in V_{2}$, then $B$ is safe
in $G$ if either
\begin{romaninline}
\item $\|B,G-B\|=0$,
\item $3\leq|B|\leq k$,
\item $|B|\geq5$,
\item $|B|=1$ and $|G|\bmod k\notin\{0,-1\}$,
\item $|B|=2$ and $|G|\bmod k\notin\{0,1\}$, or
\item $|B|=4$, $k=3$ and $|G|\bmod k\ne0$. 
\end{romaninline}
\end{cor}

\begin{proof}
As $r\in V_{2}$, $B$ is a path $v_{1}\dots v_{t}$ or a cycle $v_{1}\dots v_{t}v_{1}$
with $\|B,G-B\|=0$. If (ii) holds or both (i) and $|B|\leq k$ hold, then $v_{2}$ is a hidden
vertex, so use Corollary~\ref{cor:safe}. If (iii) holds or both (i) and $|B|\geq4$ hold, then
set $B_{0}:=v_{1}\dots v_{t-3}$ and $B_{1}:=v_{t-2}\dots v_{t}$. Now $B_{0}$ is
safe in $G_{1}:=G-B_{1}$ by Corollary~\ref{cor:d(v)=00003D1}, and $B_{1}$ is
safe in $G$ by (ii). If (iv) holds then let $f$ be an SE $L$-coloring of $G-v_{1}$
and color $v_{1}$ with $\alpha\notin T(f)\cup f(N(v_{1}))$. If (v) holds then let $f$
be an SE $L$-coloring of $G-B$. Now $|T(f)\cup f(N(v_{i}))|\leq k-3+1$ for $i\in[2]$,
so we can greedily color $v_{1}$ and $v_{2}$. If (vi) holds then $t=4$. Set $B_{1}:=v_{2}v_{3}v_{4}$.
By Corollary~\ref{cor:d(v)=00003D1}, $B_{0}:=\{v_{1}\}$ is safe in $G_{1}:=G-B_{1}$,
and by (ii) $B_{1}$ is safe in $G$. Thus $B$ is safe in $G$.
\end{proof}

\section{forbidden bugs in minimum counterexamples\protect\label{s3} }

Suppose Theorem~\ref{thm:Main4} fails for some $k\in\{3,4\}$. Then there is a
graph $G$ with $\rho_{G}^{k}\leq2-\sigma_{G}^{k}$ and a list assignment $L\in\mathcal{L}_{k}$
such that $G$ does not have an SE $L$-coloring. Choose such a counterexample $G=(V,E)$
and $L$ so that $|G|+|V_{1}(G)|$ is minimum. We emphasize that $G$ and $L$ now
form a minimum counterexample. In the rest of the proof we will obtain a contradiction
by using the minimality of $G$ to construct an SE $L$-coloring of $G$. 

Let $B\subseteq G$ be an induced subgraph, and set $G_{0}:=G-B$ .  Call $G_{0}$ \emph{handy} if $V_{1}(G_{0})\subseteq V_{1}(G)$.
The minimality of $G$ implies that $G_{0}$ is SE $L$-colorable if it is handy.
\begin{lem}
\label{lem:simple} 
\begin{enumerate}[label=(\alph{enumi})]
\item No component $B$ of $G$ satisfies $\Delta(B)\leq2$.
\item $G$ has at most one leaf; so $\sigma_{G}^{4}=|V_{1}|$.
\item If $B\subset G$ is safe in $G$, then $G_{0}:=G-B$ is not handy. 
\end{enumerate}
\end{lem}

\begin{proof}
(a) Suppose $B\subset G$ is a component of $G$ with $\Delta(B)\leq2$. Then $B$
is a cycle or path, so $B$ has $0$ or 2 leaves. 
Thus $\sigma_{G_{0}}^{k}=\sigma_{G}^{k}$. 
Since $B$ is a component of $G$,
both $G_{0}:=G-B$ is handy and  $\rho_{G_0}(A)=\rho_G(A)$ for all $A\subseteq V(G_0)$. Now we have $\rho_{G_{0}}^k\leq\rho_{G}^k\leq2-\sigma_{G}^{k}\leq2-\sigma_{G_{0}}^{k}$, and thus  $G_{0}$
has an SE $L$-coloring. By Corollary~\ref{cor:d(v)=00003D2}(i), $B$ is safe in
$G$, so $G$ has an SE $L$-coloring, a contradiction.

(b) Suppose $G$ has distinct leaves $x,y$ with $xx',yy'\in E$. By (a), $xy\notin E$.
Set $G^{+}:=G+xy$. Now $V_{1}(G^{+})=V_{1}(G)\smallsetminus\{x,y\}$, so $\sigma_{G^{+}}^{k}=\sigma_{G}^{k}$.
Let $A\subseteq V$ with $\rho_{G^{+}}^{k}=\rho_{G^{+}}^{k}(A)$. Now
\[
2-\sigma_{G}^{k}\geq\rho_{G}^{k}\geq\rho_{G}^{k}(A)\geq\rho_{G^{+}}^{k}(A)+\frac{\varepsilon_{k}}{2}|A\cap\{x,y\}|-\varepsilon_{k}\|G^{+}[\{x,y\}]\|=\rho_{G^{+}}^{k}(A)=\rho_{G^{+}}^{k}.
\]
As $|G^{+}|+|V_{1}(G^{+})|<|G|+|V_{1}(G)|$, minimality implies $G^{+}$, and so
$G$ has an SE $L$-coloring, a contradiction. 

(c) If $G_{0}$ is handy then $V_{1}(G_{0})\subseteq V_{1}(G)$, so $\sigma_{G_{0}}^{k}=_{\text{(b)}}|V_{1}(G_{0})|\leq|V_{1}(G)|=_{\text{(b)}}\sigma_{G}^{k}$
and $\rho^k_{G_{0}}(A)\leq\rho^k_{G}(A)$ for all $A\subseteq V(G_{0})$. Thus $\rho_{G_{0}}^k\leq\rho_{G}^k\leq2-\sigma_{G}^{k}\leq2-\sigma_{G_{0}}^{k}$.
By minimality, $G_{0}$ has an SE $L$-coloring $f$. As $B$ is safe in $G$, there
is an extension of $f$ to an SE $L$-coloring of $G$, a contradiction.
\end{proof}

 Let $\lambda:=\lambda_{B}$
be the number of plain $2$-threads $P\subseteq N[B]$ containing $r$, and let $\pi=\pi_{B}$
be the number of body vertices in $B$.
A maximal bug $B$ with root $r\in V_{k}$, $d_{B}(r)=d_{G}(r)$, $|B|=k+1$, and $\lambda=0$
is a \textit{wishbone} if $\pi=2$ and a \emph{jellyfish }if $\pi=1$.
See Figure~\ref{fig:Bug} (a) for a wishbone and (b) for a jellyfish when $k=4$. Each of them has $5$ vertices and no plain $2$-thread containing $r$. 
\begin{lem}
If $k$=3 then $G$ has no wishbone.\label{lem:wish}
\end{lem}

\begin{proof}
Suppose $B\subseteq G$ is a wishbone with root $r$ and leg vertex $l$. Set $B':=B-l$.
Note that the body vertices in $B$ must be in a closed thread since $G$ has at most one vertex of degree $1$.
By Lemma~\ref{lem:P+}, $B'$ is safe. For a contradiction, we will show that $G_{0}:=G-B'$
has an SE $L$-coloring. Now $V_{1}(G_{0})\smallsetminus V_{1}(G)=\{l\}$. Pick $A\subseteq G_{0}$
with $\rho_{G_{0}}^{3}=\rho_{G_{0}}^{3}(A)$. If $l\notin A$ then $\rho_{G_{0}}^{3}(A)=\rho_{G}^{3}(A)$;
else $l\in A$ and
\[
\rho_{G_{0}}^{3}(A)=\rho_{G}^{3}(A\cup B')+\frac{\varepsilon_{3}}{2}+\nu_{3}|B'|-\varepsilon_{3}(\|B'\|+\|l,B'\|)=\rho_{G}^{3}(A\cup B')+3\nu_{3}-3.5\varepsilon_{3}=\rho_{G}^{3}(A\cup B'),
\]
so $\rho_{G_{0}}^{3}=\rho_{G_{0}}^{3}(A)=\rho_{G}^{3}(A\cup B')\leq\rho_{G}^{3}\leq2-\sigma_{G}^{3}=2=2-\sigma_{G_{0}}^{3}$.
By the minimality of $G$, there is an SE $L$-coloring of $G_{0}$, a contradiction.
\end{proof}
\begin{lem}
\label{lem:jellyfish}If $k=3$ then $G$ has no jellyfish.
\end{lem}

\begin{proof}
Suppose $B\subseteq G$ is a jelly fish with root $r$, body vertex $b$ and leg
vertices $l_{1}$ and $l_{2}$. For $i\in[2]$, define $B_{i}:=B-l_{i}$ and $G_{i}:=G-B_{i}$.
Then ({*})~$V_{1}(G_{i})\smallsetminus V_{1}(G)=\{l_{i}\}$. Now each $B_{i}$ is
safe in $G$ by Corollary~\ref{cor:safe}. As $G$ has no SE $L$-coloring, for
both $i\in[2],$ $G_{i}$ has no SE $L$-coloring, and so $\rho_{G_{i}}^{3}\geq3$
by the minimality of $G$. Pick $A_{i}$ with $\rho_{G_{i}}^{3}(A_{i})=\rho_{G_{i}}^{3}$
for $i\in[2]$. As $\rho_{G_{i}}^{3}(A_{i})>\rho_{G}^{3}(A_{i})$, $l_{i}\in A_{i}$.
By supermodularity, 
\[
0\leq\rho_{G_{1}}^{3}(A_{1})+\rho_{G_{2}}^{3}(A_{2})-6=_{\text{(*)}}\rho_{G}^{3}(A_{1})+\rho_{G}^{3}(A_{2})\leq\rho_{G}^{3}(A_{1}\cup A_{2})+\rho_{G}^{3}(A_{1}\cap A_{2}),
\]
so $-\rho_{G}^{3}(A_{1}\cap A_{2})\leq\rho_{G}^{3}(A_{1}\cup A_{2})$. Using $\rho_{G}^{3}(A_{1}\cap A_{2})\leq2$,
this yields the contradiction,
\[
\rho_{G}^{3}\geq\rho_{G}^{3}(A_{1}\cup A_{2}\cup B)=\rho_{G}^{3}(A_{1}\cup A_{2})+6\|B\|-7|\{r\}|-4|\{b\}|\geq7-\rho_{G}^{3}(A_{1}\cap A_{2})\geq5.\qedhere
\]
\end{proof}
\begin{lem}
\label{lem:W=000026J-4}If $k=4$ then $G$ has no wishbone and no jellyfish.
\end{lem}

\begin{proof}
Suppose $r\in V_{4}$, $B:=B(r)$, $|B|=5$, $\lambda=0$ and $N(r)=[v_{1},\dots,v_{4}]$,
where $v_{1},v_{2}$ are leg vertices, $v_{4}$ is a body vertex and $v_{3}$ may
be either; so $B$ is a wishbone or jellyfish, depending on whether $v_{3}$ is a
body vertex. Set $P:=B-v_{1}-v_{2}$ ($P$ is a cycle or a path) and $G':=G-P+v_{1}v_{2}$.
Now $d_{G'}(x)=d_{G}(x)$ for all $x\in V(G')\smallsetminus N(v_{3})$. As $\lambda=0$,
any such neighbor $a$ satisfies $d_{G'}(a)\geq d_{G}(a)-1\geq2$. Anyway, $\sigma_{G'}^{4}\leq\sigma_{G}^{4}$.
Pick $A\subseteq G'$ with $\rho_{G'}^{4}=\rho_{G'}^{4}(A)$. If $v_{1}v_{2}\notin E(A)$
then $\rho_{G'}^{4}(A)=\rho_{G}^{4}(A)$. Else $v_{1}v_{2}\in E(A)$ and  
\[
\rho_{G'}^{4}(A)=\rho_{G}^{4}(A\cup P)+5|P|-2|V_{1}(P)|-4(\|P\|+\|P,A\|-|\{v_{1}v_{2}\}|)\leq\rho^{4}(A\cup P)\leq\rho_{G}^{4}.
\]
Now $|P|=3$ and $\|P\|+\|P,A\|=5$ regardless of whether $\pi=1$ or $\pi=2$ since
$|V_{1}(G)|\leq1$. Thus $\rho_{G'}^{4}(A)\leq\rho_{G}^{4}(A\cup P)-1\leq2-\sigma_{G}^{4}\leq2-\sigma_{G'}^{4}.$
By the minimality of $G$, there is an SE $L$-coloring $f'$ of $G'$. 
As $P$ is safe in $G$ by Corollary~\ref{cor:safe}, there is an SE
$L$-coloring of $G$, a contradiction.
\end{proof}
Call a nonempty subset $A\subseteq V$ \emph{extreme} if $\rho_{G}(A)=\rho_{G}$.
Suppose $A$ and $B$ are extreme. Using $\rho_{G}^{k}\geq\rho_{G}^{k}(A\cup B),\rho_{G}^{k}\geq\rho_{G}^{k}(A\cap B)$
and supermodularity, yields that $A\cup B$ is extreme:
\[
\rho_{G}^{k}\geq\rho_{G}^{k}(A\cup B)\geq\rho_{G}^{k}(A)+\rho_{G}^{k}(B)-\rho_{G}^{k}(A\cap B)\geq2\rho_{G}^{k}-\rho_{G}^{k}=\rho_{G}^{k}.
\]
Define $X:=\bigcup\{A\subseteq V:A\text{ is extreme}\}.$ Now $X$ is extreme and
contains every extreme set.

A \emph{fork }is a maximal bug  $B:=B(r)$ with root $r\in V_{k-1}\cup V_{k}$, $|B|=d(r)+2$,
$\lambda=1$ and $\pi=0$. In this case, $B$ is a $d(r)$-fork. 
See  Figure~\ref{fig:Bug}(c) for a $3$-fork with $k=4$. 
{When $k=3$, a $2$-fork is a $P_4$, where $r$ is an inner vertex.}

{Recall that for a path $P=v_{1}\dots v_{s}$, $\mathring{P}:=v_{2}\dots v_{s-1}$.}

{\color{black} If $v\in B$,
$w\in G_{0}:=G-B$ and $vw\in E$, then $w$ is a \emph{boarder} vertex of $G_{0}$
and an \emph{anchor} of $v$. Anchors of vertices in $B$ are called \emph{anchors}
of $B$. Suppose $A\subseteq V(G_{0})$ with $\rho_{G_{0}}(A)>\rho_{G}(A)$. Then
there is a leaf $l\in V_{1}(G_{0})\smallsetminus V_{1}(G)$; so $l$ is a boarder
vertex of $G_{0}$.}
\begin{lem}
\label{lem:4-thread}
{Let $k=3$. 
If $B$ is a $2$-fork with root $r$, 
then $N[B]\subseteq X$
and $\rho_{G}^{3}=2$.}
\end{lem}

\begin{proof}
 {\color{black} Let $B'=a_{1}v_{1}\dots v_{4}a_{2}:=N[B]$, where $r=v_{2}$ and $a_{1},a_{2}$ are (possibly equal) anchors of $B$.} 
Here we use that forks are maximal bugs. Set $R:=v_{2}v_{3}v_{4}$ and $H:=G-R$.
Now $V_{1}(H)=V_{1}(G)+v_{1}$. As $R$ is safe by Corollary~\ref{cor:safe}, $H$ has no SE $L$-coloring. By
the minimality of $G$, $\rho_{H}^{3}\geq3>\rho_{G}^{3}$. Pick $A\subseteq H$ with
$\rho_{H}^{3}=\rho_{H}^{3}(A)>\rho_G(A)$. Then $v_1\in A$. {Now $a_1\in A$, since otherwise $\rho_H(A-v_1)>\rho_H(A)$.}
Thus $a_{1}v_{1}\subseteq A$, and $\rho_{G}^{3}(A)=\rho_{H}^{3}(A)-3$.
Also set $R':=v_{1}v_{2}v_{3}$ and $H':=G-R'$. By symmetry, $V_{1}(H')=V_{1}(G)+v_{4}$,
and there is $A'\subseteq H'$ with $\rho_{H'}^{3}(A')\geq3$, $a_{2}v_{4}\subseteq A'$
and $\rho_{G}^{3}(A')=\rho_{H'}^{3}(A')-3$. By supermodularity, 
\[
0\leq\rho_{H}^{3}(A)+\rho_{H'}^{3}(A')-6=\rho_{G}^{3}(A)+\rho_{G}^{3}(A')\leq\rho_{G}^{3}(A\cup A')+\rho_{G}^{3}(A\cap A'),
\]
so $-\rho_{G}^{3}(A\cap A')\leq\rho_{G}^{3}(A\cup A')$. Now {\color{black}$\rho_{G}^{3}(A\cap A')\leq2$, so $-\rho_{G}^{3}(A\cap A')\geq-2$. Thus 
\begin{align*}
2\ge\rho_{G}^{3} & \geq\rho_{G}^{3}(A\cup A'\cup B)=\rho_{G}^{3}(A\cup A')+6\|B\|-7|\{v_2,v_3\}|\geq -\rho_{G}^{3}(A\cap A')+4\geq2.
\end{align*}
Hence $\rho_{G}^{3}(A\cup A'\cup B)=2$ and $N[B]\subseteq A\cup A'\cup B\subseteq X$.}
\end{proof}

{If $v\in B$,
$w\in G_{0}:G-B$ and $vw\in E$, then $w$ is a \emph{boarder} vertex of $G_{0}$
and an \emph{anchor} of $v$. Anchors of vertices in $B$ are called \emph{anchors}
of $B$. Suppose $A\subseteq V(G_{0})$ with $\rho_{G_{0}}(A)>\rho_{G}(A)$. Then
there is a leaf $l\in V_{1}(G_{0})\smallsetminus V_{1}(G)$; so $l$ is a boarder
vertex of $G_{0}$.}
\begin{lem}
\label{lem:Fork3,4}If $B\subseteq G$ is a fork with root $r\in V_{3}\cup V_{4}$,
then $r\in X$ and $\rho_{G}^{k}=2-\sigma_{G}^{k}$.
\end{lem}

\begin{proof}
Let the legs of $B$ be $rv_{1},\dots,rv_{d(r)-1},rxy$. If $a\in V_{3^{+}}$
is the anchor of $y$, and if $d(a)<d(r)$, then pick notation so that $av_{d(r)-1}\notin E$.
By Lemma~\ref{lem:P+}, $B':=B-v_{1}-v_{2}$ is safe in $G$. Thus $G_{0}:=G-B'$
has no SE $L$-coloring, and so $G':=G_{0}+v_{1}v_{2}$ has no SE $L$-coloring.
Using $d_{G}(a)\geq3$ and the choice of notation, $|N(a)\smallsetminus B'|\geq2$,
so $V_{1}(G')=V_{1}(G)$ and $\sigma_{G'}^{k}=\sigma_{G}^{k}$. By the minimality
of $G$, 
\begin{equation}
\rho_{G'}^{k}-1\geq2-\sigma_{G'}^{k}=2-\sigma_{G}^{k}\geq\rho_{G}^{k}.\label{eq:g'}
\end{equation}
Pick $A\subseteq V(G')$ with $\rho_{G'}^{k}(A)=\rho_{G'}^{k}$. Now $v_{1}v_{2}\subseteq A$:
else $\rho_{G'}^{k}=\rho_{G'}^{k}(A)=\rho_{G}^{k}(A)\leq\rho_{G}^{k}$, contradicting
(\ref{eq:g'}). Thus
\[
\rho_{G}^{k}\geq\rho_{G}^{k}(A+r)=\rho_{G}^{k}(A)+2\varepsilon_{k}-\nu_{k}=\rho_{G'}^{k}(A)+\varepsilon_{k}-\nu_{k}=\rho_{G'}^{k}-1\geq\rho_{G}^{k}.
\]
So $\rho_{G}^{k}=\rho_{G}^{k}(A+r)=\rho_{G'}^{k}-1$. Thus by (\ref{eq:g'}), $\rho_{G}^{k}=2-\sigma_{G}^{k}$
and $r\in X$.
\end{proof}
The proof of the next lemma involves a special case of issues we address in Section~4.
\begin{lem}
\label{lem:3+thread}Let $P\subseteq G$ be a thread with anchor $r\in V_{3^{+}}$
and $P-V_{3^{+}}=x_{1}\dots x_{t}$.
\begin{enumerate}[label=(\alph{enumi}),nosep]
\item If $P$ is \textup{loose then $t=1$ and $|G|\bmod k=0$;}
\item if $P$ is plain then 
\begin{romaninline}
\item $t=1$ and $|G|\bmod k\in\{-1,0\}$ or 
\item $t=2$ and $|G|\bmod k\in\{0,1\}$ or 
\item $t=4$, $k=3$, $P\subseteq X$, $|G|\equiv0\mod k$ and $\rho_{G}^{3}=2$; and
\end{romaninline}
\item if $P$ is closed then $d(r)\geq4$ and either
\begin{romaninline}
\item $t=2$ and $|G|\bmod k\in\{0,1\}$ or
\item $k=3$, $t=4$, $|G|\equiv0\mod k$, $P\subseteq X$ and $\rho_{G}^{3}=2$.
\end{romaninline}
\end{enumerate}
\end{lem}

\begin{proof}
Suppose not. We will construct an induced subgraph $S$ such that $S$ is safe in $G$ and $G_{0}:=G-S$
is handy, thus contradicting Lemma~\ref{lem:simple}(c). 

If (a) fails then $P$ is loose and either $t\geq2$ or both $t=1$ and $|G|\bmod k\ne0$.
Set {\color{black}$S:=P-r$}. By Corollary~\ref{cor:d(v)=00003D1}, $S$ is safe in $G$;
and $G_{0}$ is handy since $d_{G_{0}}(r)=d_{G}(r)-1\geq2$. 

If (b) fails then $P$ is plain and all of (bi--biii) fail. If (biii) fails due to $k=4$, then $t\neq 4$ otherwise $P$ by Lemma~\ref{lem:P+} and $G_0$ is handy by the definition of closed thread; thus $k=3$, and Lemma~\ref{lem:4-thread}
implies $t\ne4$. {\color{black}Set $S:=P-\{r,a\}$, where  $a$ is the other anchor of $P$.} As (bi--ii) fail, $S$ is safe in $G$
by Corollary~\ref{cor:d(v)=00003D2}(ii,iii). Also, $G_{0}$ is handy since $d_{G_{0}}(r)=d_{G}(r)-1\geq2$
and $d_{G_{0}}(a)=d_{G}(a)-1\geq2$, where $a\ne r$ is the other end of $P$. 

If (c) fails then $P$ is closed. Suppose $d(r)=3$. Set $S:=B(r)$. If $S=P\cup R$ for some loose thread $R$, then $S=G$ and by (a) $k=3$ and $|S|=|G|=4$.
Clearly $G-P$ has an SE $L$-coloring and $P$ is safe in $G$ by Lemma~\ref{lem:P+}.
Now it must be the case that $S=P\cup R$
for some $r,a$-plain-thread $R$.
Let $r'\in R\cap N(r)$, so $R$ is a maximal bug rooted at $r'$. As $\|S,G_{0}\|\leq1$, $G_{0}$ is handy. Thus
$S$ is not safe. By Corollary~\ref{cor:d(v)=00003D2} applied to $r'$ and $R$, either (ci) holds or $t = 4$ and $k = 3$.
If (ci) then by Corollary~\ref{cor:safe}, $P$ is safe in $G$ and $|S|\geq k+1$.
By Lemmas~\ref{lem:wish} and \ref{lem:W=000026J-4}, $S$ is not a wishbone, so
$|S|\geq k+2$. Thus {\color{black}$|R-r|\geq2$.} By Corollary~\ref{cor:d(v)=00003D1},
{\color{black}$R-r$} is safe in $G-P$, so $S$ is safe in $G$, a contradiction. Else
$t=4$ and $k=3$. By Corollary~\ref{cor:d(v)=00003D1}, {\color{black}$R+rx_{4}$} is
safe in $G-x_{1}x_{2}x_{3}$, and by Corollary~\ref{cor:safe}, $x_{1}x_{2}x_{3}$
is safe in $G$, a contradiction. So $d(r)\geq4$.

Set $S:=P-r$. As $\|S,G_{0}\|\leq2$ and $d(r)\geq4$, $G_{0}$ is handy. If $S$
is safe then this is a contradiction. Thus $S$ is not safe. As $P$ is a maximal bug for some $r'\in P\cap V_2$, by Corollary~\ref{cor:d(v)=00003D2}(ii,iii), $t=2$ holds or $t=4,k=3$ and $|G|\equiv0\bmod 3$ hold.
For the latter, Lemma~\ref{lem:4-thread},
$P\subseteq X$ and $\rho_{G}^3=2$, so (cii) holds.
If $t=2$, then either $S$ is safe by corollary 4.5(v) or (ci) is true.
\end{proof}

\section{Special forbidden bugs in minimum counterexamples}

Let $B\subseteq G$ be an induced subgraph, and set $G_{0}:=G-B$.
{Recall that $v\in G_0$ is a boarder vertex for $G_0$ if it has a neighbor $w\in B$, and that in this case, $v$ is an anchor of $w$.}
{An \emph{anchor
for a vertex of $B$ is a boarder vertex} for $G_{0}$.}
By Lemma~\ref{lem:simple}(c),
if $B$ is safe and $G_{0}$ is handy, then $G$ has an SE $L$-coloring. But if $B$
has anchors that are leaves in $G_{0}$, then $G_{0}$ is not handy. In this section
we deal with this situation. 

\subsection{Bugs, buffers and cores }

 Suppose $B\subseteq G$ and $G_{0}:=G-B$. By Lemma~\ref{lem:simple}(b), $|V_{1}(G)|\leq1$.
If $V_{1}(G)\ne\emptyset$ then let $l$ be the unique leaf of $G$; else (for notational
purposes) let $l\notin G$. For each component $Z$ of $G_{0}$, pick $H_{Z}\subseteq Z$
as follows. If $Z$ is a tree then set $H_{Z}=\{y\}$ for some vertex $y\in Z$,
preferring $y=l$. {\color{black}Else $Z$ is not a tree; as $Z$ is a component of $G_0$, $Z$ has a cycle.}  Let $H'_{Z}$ be the union of all $Y\subseteq Z$ with $\delta(Y)\geq2$.
{\color{black}If $l$ is in $Z$, but $l\notin H'_{Z}$, then let $P$ be an
$l,H'_{Z}$-path,} and set $H_{Z}:=H_{Z}'\cup P$. Finally, put $H:=\bigcup H_{Z}$.
Now for all $v\in V(H)$, either $d_{H}(v)=0$ or {\color{black}$d_{H}(v)\geq2$ or $v$ is the unique leaf $l$ of $G$. Thus}
$H$ is handy, $\sigma_H = \sigma_G$ and $\rho_H \leq \rho_G$. {\color{black}We call $H$ the \emph {core} of $G_0$.} 

Set $F:=(G_{0}-H)+N(G_{0}-H)\cap H+E(G_{0}-H,H)$. So {\color{black}$G_{0}= F\cup H$.}
Now $F$ is acyclic since any cycle in $F$ would be added to $H$. Also every component
of $F$ has exactly one vertex in $H$: at most one since any path in $F$ between
two vertices of $H$ would be added to $H$; at least one since every component of
{\color{black}$G_0$} meets $H$. {We call $F$ the \emph{buffer} of $G_0$ between $B$ and $H$}.

\begin{figure}
\begin{tikzpicture}[scale=1.2]

\draw[rounded corners, color=green, fill=green!5,fill opacity=0.6](-.7,2.7) rectangle (1.1,.9);

\draw[rounded corners, color=blue, fill=blue!5,fill opacity=0.6](-1.2,1.2) rectangle (3.2,-.1);

\def \vt{circle (1.5pt) [fill]} 

\draw (1,4)--(.5,4.5);
\draw (1,4)--(1.5,4.5);
\draw(1.5,4.5)--(.5,4.5);
\draw (1,4) \vt;
\draw (1.5,4.5) \vt;
\draw (.5,4.5) \vt;
\draw (.7,4) node(r) {$r$};

\draw (1,4)--(-1,3);
\draw (-1,3) \vt;
\draw (0,3.5) \vt;

\draw (2.5,3) \vt;
\draw (1,4)--(2.5,3)--(2.5,1);
\draw (2.5,1) \vt;
\draw (3,.5) \vt;
\draw (2,.5) \vt;
\draw (2.5,0) \vt;
\draw (2,.5)--(2.5,1)--(3,.5)--(2,.5)--(2.5,0)--(3,.5);

\draw (-.5,1) \vt;
\draw (-1,.5) \vt;
\draw (0,.5) \vt;
\draw (-.5,0) \vt;
\draw (0,.5)--(-.5,1)--(-1,.5)--(0,.5)--(-.5,0)--(-1,.5);

\draw (.3,3) \vt;

\draw (-.25,2.25) \vt;
\draw (1,4)--(.3,3)--(-.25,2.25)--(-.5,1.5);
\draw (.65,3) \vt;
\draw (1,4)--(.65, 3)--(-.25,2.25);

\draw (-.5,2.5) \vt;
\draw (0,3) \vt;
\draw (-1,3)--(-.5,2.5);
\draw (1,4)--(-.5,2.5);
\draw (-.5,2) \vt;
\draw (-.5,1.5) \vt;
\draw (-.5,1)--(-.5,1.5)--(-.5,2.5);

\draw (2.2,3)\vt;
\draw (1.8,3)\vt;
\draw (1.4, 3) \vt;
\draw (1,4)--(1,.7);
\draw (1,.7) \vt;
\draw (.7,.7) node(v) {$v$};
\draw (1,2) \vt;
\draw (1,4)--(1.4,3)--(1,2);
\draw (1,4)--(1.8,3)--(1,.7);
\draw (1,4)--(2.2,3)--(1,.7);

\draw (0,0) \vt;
\draw (.5,0) \vt;
\draw (-.5,0)--(0,0)--(.5,0);
\draw (.5,.3) node(l) {$\ell$};
\draw[thick,decorate,decoration={brace,amplitude=8pt,mirror}] (-.5,-.1) -- (.5,-.1) node[midway,below=8pt] {$P$};

\draw[red,dashed,thick] (-2.0,2.8) -- (3.6,2.8);

\node[red!60!black] at (-2,4) {\large\textbf{B}};        
\node[green!60!black] at (-2,2.2) {\large\textbf{F}}; 
\node[blue!60!black] at (-2,0.6) {\large\textbf{H}};       

\draw[thick,decorate,decoration={brace,amplitude=8pt,mirror}] (3.8,-.1) -- (3.8,2.7) node[midway,right=10pt] {$G_0$};

\end{tikzpicture}
\caption{\label{fig:Buffer} {The maximal bug $B(r)$ is above the dotted line; the buffer $F$ is a forest with two components; the core $H$ has three components,  one of which contains the path $P$ ending in the only leaf $\ell$ of $G$.}} 
\end{figure}

\begin{lem}
\label{lem:buf}$Let$ $G_{0}$ be a proper induced subgraph of $G$ with core $H$
and buffer $F$. Then there is an independent set $Q\subseteq V_{1}(F)\smallsetminus H$
such that  $G_{0}-Q$ has an SE $L$-coloring.
\end{lem}

\begin{proof}
As $H$ is handy, $H$ has an SE $L$-coloring $f$. Arguing by induction on $|F|$
we show that there is a set $Q\subseteq V_{1}(F)\smallsetminus H$ such that $f$
can be extended to an SE $L$-coloring of $G_{0}-Q$. The case $|F|=0$ is trivial,
so suppose $|F|\geq1$. Let $x\in V_{1}(G_{0})\smallsetminus H$ be a leaf of $F$.
Set $B':=B_{F}(x)-H$, $F':=F-B'$ and $G':=G_{0}-B'$. By induction there is a set
$Q'\subseteq V_{1}(F')\smallsetminus H$, such that $f$ can be extended to an $SE$
$L$-coloring $f'$ of $G'-Q'$. If $|B'|=1$ then the anchor $a$ of $x$ in $G'$
is colored by $f'$, and so is not in $Q'$. Thus $Q:=Q'+x$ is an independent subset
of $V_{1}(F)-H$. Else set $Q:=Q'$; by Corollary~\ref{cor:d(v)=00003D1}, $f'$
can be extended to $(G'-Q')\cup B'=G_{0}-Q$.
\end{proof}
Suppose $B$ is a maximal bug in $G$. By Lemma~\ref{lem:buf}, there is an independent
subset $Q\subseteq N(B)\cap V_{1}(G-B)$ such that $G':=G-B-Q$ has an SE $L$-coloring
$f$. We call $(G',B,Q,f)$ a \emph{good string} in $G$; note that $Q$ could be
empty. As $B$ is a maximal bug, every vertex $q\in Q$ satisfies $d(q)\geq3$ and
$\|q,B\|\geq2$. 

\subsection{Applications}

Fix $r\in V_{3^{+}}\smallsetminus X$ such that $r$ is not the root of a fork. Define $B:=B(r)$ and
$G_{0}:=G-B$. Let $(G',B,Q,f)$ be a good string in $G$ (maybe $Q=\emptyset$).
Now $\|B,G_{0}\|=d(r)-\pi$. Every $q\in Q$ satisfies $\|q,B\|\geq2$, so $|Q|\leq(d(r)-\pi)/2$;
if $Q\ne\emptyset$ then let $q\in Q$ and let $a$ be the anchor of $q$ in $G'$. 

By Lemma~\ref{lem:3+thread}, if $|B(r)|\geq 2$, then there are no $3^{+}$-threads or $2^+$-loose threads incident to $r$, and
so 
\begin{equation}
|B(r)|\leq2d(r)+1-\pi.\label{eq:bigdegree}
\end{equation}

\begin{lem}
\label{prop:k=00003D3=00003Dr}If $k=3=d(r)$ and $r$ is not the root of a fork,
then $|B|\leq4-\pi$.
\end{lem}

\begin{proof}
Suppose $|B|\geq5-\pi$. If $|B|=3$ then $\pi=2$ and $\|B,G_{0}\|=1$. Thus $G_{0}$
is handy and $B$ is safe by Corollary~\ref{cor:safe}, a contradiction towards Lemma~\ref{lem:simple}(c). If $|B|=4$
then $\pi\geq1$. By Lemmas~\ref{lem:jellyfish} and \ref{lem:wish}, $B$ is neither
a jellyfish nor wishbone. Thus $B$ is a 4-loose-thread, contradicting Lemma~\ref{lem:3+thread}(a). 

Suppose $|B|=5$. Now $\lambda\geq1$. By assumption, $B$ is not
a fork, so $\lambda\geq2$ or $\pi\geq1$. Now $B$ is safe by Corollary~\ref{cor:bsafe}.
Thus $G_{0}$ is not SE $L$-colorable, so as $|Q|\le (3-\pi)/2\le 1$, we may set $Q=\{q\}$. Let $B_{0}=yx$,
where $P=aqyxr$ and $yxr\subseteq B$. Set $B_{1}:=(B-B_{0})\cup Q$. Using Corollary~\ref{cor:safe},
$B_{0}$ is safe in $G_{1}:=G-B_{1}-Q$ and $B_{1}$ is safe in $G$, a contradiction. 

Suppose $|B|=6$. Now $\lambda=2$. Say $V(B)=\{r,x_{1},y_{1},x_{2},y_{2},x\}$ where
$rx_{1}y_{1}$ and $rx_{2}y_{2}$ are legs. By Corollary~\ref{cor:bsafe}, $B$
is safe, so $Q=:\{q\}$. If $qx\in E$ then set $B_{0}:=qx$ and $B_{1}:=B-B_{0}$.
By Corollary~\ref{cor:safe}, the 2-loose-thread $aqx$ is safe in $G_{1}:=G-B_{0}$,
and by Corollary~\ref{cor:bsafe}, $B_{1}$ is safe in $G$. Else $qy_{1},qy_{2}\in G$.
Set $B_{0}:=xrx_{1}y_{1}$ and $B_{1}:=(B-B_{0})\cup Q$. Then $B_{0}$ is the interior of
a $4$-loose-thread or a $4$-closed-thread, it is safe in $G_{1}:=G-B_1$, and $B_{1}$ is safe in $G$
by Corollary~\ref{cor:safe}. Anyway, $f$ can be extended to an SE $L$-coloring
of $G$, a contradiction.

Suppose $|B|=7$. Now $\lambda=3$. Set $B_{0}:=Q\cup y_{1}x_{1}$, where $y_{1}x_{1}r$
is a leg of $B$, and define $B_{1}:=B-B_{0}$. Now $B_{0}$ is safe in $G_{1}:=G-B_{1}$
by Corollary~\ref{cor:safe}, and $B_{1}$ is safe in $G$ by Corollary~\ref{cor:bsafe}.
This is a contradiction.
\end{proof}
\begin{lem}
\label{lem:k=00003D4,d(v)=00003D3}If $k=4$ and $d(r)=3$ and $r$ is not the root
of a fork, then $|B(r)|\leq3-\pi$.
\end{lem}

\begin{proof}
Suppose $|B|\geq4-\pi$. Using $|B|\geq1+\pi$ and (\ref{eq:bigdegree}), $3\leq|B|\leq7$.
As $\|B,G_{0}\|\leq3-\pi$, Corollary~\ref{cor:safe} implies $B$ is safe, if $|B|\leq4$.
By Lemma~\ref{lem:W=000026J-4}, $B$ is neither a wishbone nor jellyfish. As $B$
is not a fork, if $|B|\geq5$ then $\lambda+\pi/2>1$. By Corollary~\ref{cor:bsafe},
$B$ is still safe. Anyway, $G_{0}$ is not SE $L$-colorable, so $Q=\{q\}$. For
a contradiction it suffices to show that $B^{+}:=B+q$ is safe in $G$. 

Set $G':=G-B^{+}$. Now $\|B,G'\|\leq d(r)-\|B,q\|-\pi(B)\leq1-\pi(B)$. So $\pi(B)\leq1$
and $|B^{+}|\geq4$. If $|B^{+}|=4$ then $B^{+}$ is safe in $G$ by Lemma~\ref{lem:P+}.
Else $|B^{+}|\geq5$ and there are at least three vertices in $B$ with no neighbors
in $G'$; as $\pi\leq1$, two of these vertices, say $x,y$, are nonadjacent. If
$|B^{+}|=5$ then set $B_{1}:=\{x\}$ and $B_{2}:=B^{+}-B_{1}$. Now $B_{1}$ is
safe in $G_{1}:=G-B_{2}$. As $\|B_{2},G_{1}\|\leq\|q,G'\|+\|B,G'\|+d(x)\leq4$,
$xy\notin E$, $1\leq\|q,G_{1}\|\leq2$ and $|B_{2}|=4$, using Lemma~\ref{lem:P+},
$B_{2}$ is safe in $G$. Thus $B^{+}$ is safe in $G$. 

Else $|B^{+}|\geq6$, and there is a $4$-path $P=qx_{1}x_{2}x_{3}\subseteq B^{+}$.
Define $B_{2}:=P$ and $B_{1}:=B^{+}-B_{2}$. As $\|B_{1},G'\|\leq1$ and $2\leq|B_{1}|\leq4$,
Lemma~\ref{lem:P+} implies $B_{1}$ is safe in $G_{1}:=G-B_{2}$. As $\|q,G_{1}\|\leq\|q,G'\|+\|q,B_{1}\|\leq3$,
$\|x_{1},G_{1}\|=0$ and $\|\{x_{2},x_{3}\},G_{1}\|\leq2$, Lemma~\ref{lem:P+}
implies $B_{2}$ is safe in $G$. Thus $B^{+}$ is safe in $G$.
\end{proof}
\begin{lem}
\label{lem:d(v)=00003D4}If $k=4=d(r)$ and $r$ is not the root of a fork, then
$|B|\leq5-\pi$.
\end{lem}

\begin{proof}
Suppose $|B|\geq6-\pi$. First we show that $B$ is safe in $G$. If $|B|=4$ then
$\pi=2$, so $B$ is safe in $G$ by Corollary~\ref{cor:safe}. Suppose $|B|=5$.
Then $\pi\geq1$. By Lemma~\ref{lem:W=000026J-4}, $B$ is neither a jellyfish nor
a wishbone. Thus $\lambda\geq1$. Now $B$ is safe by Corollary~\ref{cor:bsafe}.
Suppose $|B|\geq6$. The rest of this paragraph holds even when $B$ is not maximal.
Now $\lambda\geq1$. If $\pi\geq1$ then $|B|\leq8$, so $B$ is safe in $G$ by
Corollary \ref{cor:bsafe}. Else $\pi=0$. As $B$ is not a fork, $\lambda\geq2$.
If $|B|\leq8$ then $B$ is safe in $G$ by Corollary~\ref{cor:bsafe}. Else $B=9$
and $\lambda=4$. Let $B_{0}=rx_{1}x_{2}\subseteq B$ and $B_{1}:=B-B_{0}$. Now
$B_{0}$ is safe in $G_{1}:=G-B_{1}$ by Corollary~\ref{cor:safe}, and $B_{1}$
is safe in $G$ by Corollary~\ref{cor:bsafe}, so $B$ is safe in $G$, and $G_{0}$
is not SE $L$-colorable, so $1\leq|Q|\leq2$. For a contradiction, we will show
that $B^{+}:=B\cup Q$ is safe.  Pick $Z\subseteq N(Q)\cap B$ so that $|Z|\in\{1,2\}$
and $|B-Z|\ne5$. Set $B_{1}:=B-Z$ and $B_{0}:=Q\cup Z$. By Lemma~\ref{lem:P+},
$B_{0}$ is safe in $G_{1}:=G-B_{1}$. Now $B_{1}$ is a bug. If $|B_{1}|\geq6-\pi$
then $B_{1}$ is safe by the previous paragraph; else $|B_{1}|-\pi\leq4$ and $B_{1}$
is safe by Corollary~\ref{cor:safe}. Anyway $B\cup Q$ is safe. 
\end{proof}

\section{Discharging}

We are now ready to prove Theorem~\ref{thm:Main4}. Assign a charge $\mu(x)$ to
each $x\in V\cup E$ by
\begin{align*}
\mu^{k}(v) & =\begin{cases}
-\nu_{k}+\frac{\varepsilon_{k}}{2} & \text{if }v\in V_{1}\\
-\nu_{k} & \text{if }v\in V\smallsetminus V_{1}
\end{cases}\negthickspace;~~~~~~\text{ and }~~~~~~\mu^{k}(e)=\varepsilon_{k}~~\text{if }e\in E.
\end{align*}
Now 
\[
\rho_{G}^{k}(V)=\sum_{v\in V}\mu^{k}(v)+\sum_{e\in E}\mu^{k}(e).
\]
For $Y\subseteq V$, consider the following discharging rules where $\delta_{k}:=\nu_{k}-\varepsilon_{k}=1$.
\begin{enumerate}[label=(R\arabic{enumi})]
\item  The charges of vertices and edges in $G[Y]$ do not change.
\item Every edge in $E(Y,V\smallsetminus Y)$ sends $\varepsilon_{k}$ to its end in $V\smallsetminus Y$.
\item Every edge in $G-Y$ sends $\varepsilon_{k}/2$ to each of its ends.
\item Every vertex $v\in V_{3^{+}}\smallsetminus Y$ sends $\delta_{k}/2$ to each of its
leg vertices in $B(v)\cap V\smallsetminus Y$ and $\delta_{k}$ to each of its body
vertices in $B(v)\cap V\smallsetminus Y$.
\item Each $v\in V_{2}\cap N(Y)$ sends $\delta_{k}$ to its neighbor in $V\smallsetminus Y$
if it exists.
\end{enumerate}
Let $\mu_{Y}^{k}(x)$ be the final charge of each $x\in V\cup E$. As $\mu_{Y}(e)=0$
for all edges $e\in E\smallsetminus E(G[Y])$, 
\begin{equation}
\rho_{G}^{k}\geq\rho_{G}^{k}(V)=\rho_{G}^{k}(Y)+\sum_{v\in V(G)-Y}\mu_{Y}^{k}(v).\label{eq:finc'-1}
\end{equation}
\label{eq:5.1}
\begin{lem}
\label{lem:charge}
\begin{romaninline}
\item If all roots of forks are in $Y$, then $\mu_{Y}^{k}(v)\geq0$ for all $v\in V\smallsetminus Y$.
\item $X=V$.
\end{romaninline}
\end{lem}

\begin{proof}
Consider any vertex $v\in V\smallsetminus Y$. By Lemma~\ref{lem:simple}(a), the
component $H$ of $v$ satisfies $\Delta(H)\geq3$. Thus $v$ is in at least one
bug $B=B(r_{v})$ with $r_{v}\in V_{3^{+}}$, and $v$ is not isolated. As $v$ is
not in a $4^{+}$-thread or a loose $2^{+}$-thread by Lemma~\ref{lem:4-thread},
$r_{v}$ can be chosen in $N[v]$.

Suppose $v\in V_{1}$. Then $v$ is a body vertex of $r_{v}$. If $r_{v}\in Y$ then
$v$ receives $\varepsilon_{k}$ from $vr_{v}$ by (R2). So
\[
\mu_{Y}^{k}(v)=\mu^{k}(v)+\varepsilon_{k}=(-\nu_{k}+\frac{\varepsilon_{k}}{2})+\varepsilon_{k}=\frac{\varepsilon_{k}}{2}-\delta_{k}\geq1.
\]
 Else $r_{v}\in V\smallsetminus Y$, so $v$ receives $\varepsilon_{k}/2$ from $vr_{v}$
by (R3) and $\delta$ from $r_{v}$ by (R4). Thus
\[
\mu_{Y}^{k}(v)=\mu^{k}(v)+\frac{\varepsilon_{k}}{2}+\delta_{k}=(-\nu_{k}+\frac{\varepsilon_{k}}{2})+\frac{\varepsilon_{k}}{2}+(\nu_{k}-\varepsilon_{k})=0.
\]

Suppose $v\in V_{2}$. If $v\in N(y)$ with $y\in Y$ then $v$ receives $\varepsilon_{k}$
from $yv$ by (R2), $v$ receives at least $\varepsilon_{k}/2$ from its other incident
edge by (R2,3), and $v$ sends $\delta{}_{k}$ to its neighbor in $V\smallsetminus Y$
(if it exists) by (R5). So
\[
\mu_{Y}^{k}(v)\geq\mu^{k}(v)+\varepsilon_{k}+\frac{\varepsilon_{k}}{2}-\delta_{k}=-\nu_{k}+\frac{3}{2}\varepsilon_{k}-(v_{k}-\varepsilon_{k})=\frac{\varepsilon_{k}}{2}-2\delta_{k}\geq0.
\]
Else $N(v)\subseteq V\smallsetminus Y$. Now $v$ receives $\varepsilon_{k}/2$ from
each of its incident edges by (R3). If $v\in N(w),$ where $w\in V_{2}\cap N(Y)$,
then $v$ receives $\delta_{k}$ from $w$ by (R5). Else $v\notin N(N(Y))$. If $v$
is a leg vertex of $B(r_{v})$, then $v$ receives $\delta_{k}/2$ from both $r_{v}$
and the other end of its plain thread by (R4). Else $v$ is a body vertex of $B(v_{r})$,
so $v$ receives $\delta_{k}$ from $r_{v}$ by (R4). Anyway 
\[
\mu_{Y}^{k}(v)\geq\mu^{k}(v)+2\frac{\varepsilon_{k}}{2}+\delta_{k}=-\nu_{k}+\varepsilon_{k}+\nu_{k}-\varepsilon_{k}=0.
\]

Otherwise $v\in V_{3^{+}}$. By (R2,3), $v$ receives at least $\varepsilon_{k}/2$
from each incident edge. By (R4), $v$ sends $\delta_{k}/2$ to each leg vertex of
$B(v)\smallsetminus Y$ and $\delta_{k}$ to each body vertex of $B(v)\smallsetminus Y$.
Thus 
\[
\mu_{Y}^{k}(v)\geq-\mu^{k}(v)+\frac{\varepsilon_{k}}{2}d(v)-\frac{\delta_{k}}{2}(|B|-1+\pi).
\]
As $y$ is not the root of a fork, the following table gives bounds on $\mu_{Y}^{k}(v)$
for all $v\in V_{3^{+}}\smallsetminus Y$. 
\begin{center}
\begin{tabular}{|c|c|c|c|c|c|}
\hline 
$k$ & $d(v)$ & $|B|\leq$ & $-\mu^{k}(v)+\frac{\varepsilon_{k}}{2}d(v)-\frac{\delta_{k}}{2}(|B|-1+\pi)$ & $\mu_{Y}^{k}(v)\geq$ & Reference\tabularnewline
\hline 
$3$ & $3$ & $4-\pi$ & $-7+9-1.5$ & $.5$ & Lemma~\ref{prop:k=00003D3=00003Dr}\tabularnewline
\cline{2-6}
 & $\geq4$ & $2d(v)+1-\pi$ & $-7+3d(v)-d(v)$ & $2d(v)-7$ & (\ref{eq:bigdegree})\tabularnewline
\hline 
$4$ & $3$ & $3-\pi$ & $-5+6-1$ & $0$ & Lemma~\ref{lem:k=00003D4,d(v)=00003D3}\tabularnewline
\cline{2-6}
 & $4$ & $5-\pi$ & $-5+8-2$ & $1$ & Lemma~\ref{lem:d(v)=00003D4}\tabularnewline
\cline{2-6}
 & $\geq5$ & $2d(v)+1-\pi$ & $-5+d(v)$ & $d(v)-5$ & (\ref{eq:bigdegree})\tabularnewline
\hline 
\end{tabular}
\par\end{center}

\noindent So $\mu_{Y}^{k}$$(v)\geq0$, completing the proof of (i). For (ii), let
$Y=X$. By Lemmas~\ref{lem:4-thread} and (\ref{lem:Fork3,4}), the root of every
fork is in $Y$. So by \ref{eq:finc'-1}: 
\[
\rho_{G}^{k}=\rho_{G}^{k}(X)\leq\rho_{G}^{k}(X)+\sum_{v\in V(G)-X}\mu_{X}(v)=\rho_{G}^{k}(V)\leq\rho_{G}^{k}.
\]
Thus $\rho_{G}^{k}=\rho_{G}^{k}(V)$, and $V=X$ by the maximality of $X$. 
\end{proof}
\[
\]

\begin{lem}
\label{lem:cong}
\begin{romaninline}
\item If $G$ has a leaf then $k=3$ and $\rho_{G}^{3}=0$. Anyway, 
\item $\rho_{G}\equiv-|G|\mod k$ and
\item $G$ has no fork.
\end{romaninline}
\end{lem}

\begin{proof}
Using $k\in\{3,4\}$, Lemma~\ref{lem:charge} yields 
\begin{equation}
\rho_{G}^{k}=\rho_{G}^{k}(X)=\rho_{G}^{k}(V)=\varepsilon_{k}\|G\|-\nu_{k}|G|+\frac{\varepsilon_{k}}{2}|V_{1}|\equiv-|G|+\frac{\varepsilon_{k}}{2}|V_{1}|\mod k.\label{eq:mod}
\end{equation}
Suppose $G$ has a leaf. By Corollary~\ref{cor:d(v)=00003D1}, $|G|\bmod k=0$.
By Lemma~\ref{lem:simple}(b), $|V_{1}|\leq1$, so ({*})~$\rho_{G}^{k}\equiv\varepsilon_{k}/2\mod k$.
Thus $0\leq\rho_{G}^{k}\leq2-\sigma_{G}^{k}$. If $k=4$ then $\rho_{G}^{4}=1$,
contradicting ({*}); else $k=3$ and $\rho_{G}^{3}\bmod3=0$, so (i) holds. Anyway,
$\varepsilon_{k}/2\bmod k=0$ or $|V_{1}|=0$, so (ii) holds. 

Suppose $F$ is a fork in $G$. If $G$ is a $2$-fork then $G$ has a plain $4$-thread;
if $G$ is a $3^{+}$-fork then $G$ has a plain $1$-thread and a plain $2$-thread.
Regardless, by (i,ii) and Lemmas~\ref{lem:4-thread}--\ref{lem:3+thread}, 
\[
2\equiv_{(i)}2-\sigma_{G}^{k}\equiv_{\text{L\ref{lem:4-thread},\ref{lem:Fork3,4}}}\rho_{G}^{k}\equiv_{\text{(ii)}}-|G|\equiv_{\text{L\ref{lem:3+thread}}}0\mod k,
\]
a contradiction. So (iii) holds.
\end{proof}
Now we finish the proof of Theorem~\ref{thm:Main4}. By Lemma~\ref{lem:cong}(iii),
$G$ has no forks, so by Lemma~\ref{lem:charge}, $\mu_{\emptyset}^{k}(v)\geq0$
for all $v\in V$. By Lemma~\ref{lem:simple}(a), $G$ has a $3^{+}$-vertex, say
$v$. By Lemma \ref{lem:cong}(ii), $\rho_{G}^{k}\equiv-|G|\mod k$.

First suppose $\rho_{G}^{k}=2$. Now $|G|\bmod k\ne0$, so $G$ has no leaf (or loose-thread)
by Lemma~\ref{lem:cong}(i). If $k=4$ then $|G|\bmod k=2\not\equiv-1,0,1$. By
Lemma~\ref{lem:3+thread}(b,c), $G$ has no threads, so $\delta(G)\geq3$, and we
have the contradiction,
\[
\rho_{G}^{4}=-\nu_{4}|G|+\varepsilon_{4}\|G\|\geq-\nu_{4}|G|+\varepsilon_{4}\cdot\frac{3}{2}|G|\geq|G|\geq4.
\]
Else $k=3$ and $|G|\bmod k=1$. By Lemma~\ref{lem:3+thread}, $G$ has no $1$-thread,
and each $2$-vertex has a $2$-neighbor and a $3^{+}$-neighbor. If $v\in V_{3}$
then $1+2\lambda+\pi=|B(v)|\leq4-\pi$ and $\pi\ne1$. Thus $\pi=0$ and $\lambda\leq1$.
Now $\mu_{\emptyset}^{3}(v)=1$, and $v$ has two $3^{+}$-neighbors, each with charge
at least $1$, a contradiction. So $V_{3}=\emptyset$. If $d(v)\geq5$ then $\mu_{\emptyset}^{3}(v)\geq3$,
a contradiction. So $d(v)=4$ and $\mu_{\emptyset}^{3}(v)\geq1$. Thus $V=V_{2}\cup V_{4}$,
$|V_{4}|\leq\rho_{G}^{3}$, $5\leq|N[v]|\leq|G|\leq5|V_{4}|$, $|G|\bmod3=1$ and
$|V_{2}|$ is even. So $|V_{4}|=2$ and $|V_{2}|=8$. Set $V_{4}:=\{v_{1},v_{2}\}$.
For a contradiction we will find an SE $L$-coloring of $G$. As $|G|=10$, we are
allowed one $4$-class.

If possible, color some $v_{i}\in V_{4}$ and some $w_{1},w_{2},w_{3}\in N(v_{3-i})$
with some color $\alpha$, leaving $w_{4}\in N(v_{3-i})$ uncolored. Now $|\widetilde{\alpha}|=4$.
Extend this to an SE $L$-coloring $f$ of $G$ by greedily coloring in the order
$w_{4},N(w_{4})\smallsetminus N(v_{3-i}),N(v_{i})\smallsetminus N(w_{4})$ so that
$f(w_{4})\ne\alpha$ and no color is used more than twice on $N(v_{i})$; the latter
is possible because each vertex of $N(v_{i})\smallsetminus N(w_{4})$ has at least
two available colors after $N(w_{4})$ is colored. Else, for each $v_{i}\in V_{4}$,
each color $\alpha\in L(v_{i})$ appears in the lists of at most two vertices in
$N(v_{3-i})$, so neither vertex of $V_{4}$ can appear in a $4$-class. Thus it
suffices to color $G$ so that $V_{2}$ contains no $4$-class. For this, color $w_{1},w_{2}\in N(v_{1})$
and $w_{3}\in N(v_{2})$ distinctly, and then continue greedily in the order $v_{1},v_{2},\dots$.

Now suppose, $\rho_{G}\leq1$. If $k=3$ then $\mu_{\emptyset}^{3}(v)>0$ so $\rho_{G}=1$.
Thus $V_{3^{+}}\in\{V_{3},V_{4}\}$ and $|V_{4}|\leq1$. Also, every thread in $G$
is a $1$-plain-thread, so $G=K_{2,3}$ and $G$ has an SE $L$-coloring, a contradiction.
Else $k=4$. If $\rho_{G}=0$ then $\Delta(G)\leq3$, since vertices with higher
degree get positive charge. Else $\rho_{G}=1$, so $G$ has no $2$-threads and no
loose-threads; in particular $\pi=0$. So if $r\in V_{3^{+}}$ then $|B(r)|\leq1+d(r)-\pi$
and $\mu_{\emptyset}^{4}(v)\geq1.5d(v)-5$. Thus $\Delta(G)\leq4$. If $w\in V_{4}$
is adjacent to $x\in V_{3^{+}}$, then $|B(w)|\leq4$, so $\mu_{\emptyset}^{4}(w)\geq-5+8-1.5=1.5$.
Anyway, $\theta(G):=\max_{xy\in E}(d(x)+d(y))\leq6$. By \cite[Theorem~6.4]{KKY},
$G$ has an equitable $L$-coloring $f$. As $|G|\bmod4\equiv-\rho_{G}\bmod4\in\{0,-1\}$,
$f$ is an SE $L$-coloring. 

\section{Concluding remarks}

\noindent 1. The situations with colorings in $3$ and $4$ colors are somewhat different.
While there are sharpness examples for the $3$-color part of  Theorem~\ref{thm:main0} with maximum degree
$10$ and arbitrarily many vertices, we do not know sharpness examples for the $4$-color part
 with maximum degree less than half of the number of vertices.

\medskip{}

\noindent 2. We expect that the largest maximum average degree guaranteeing that
a graph with minimum degree at least $2$ is SE $5$-choosable is around $\frac{14}{5}$.
But we do not have a good guess what is the largest maximum average degree that provides
that a graph with minimum degree at least $2$ is SE $6$-choosable or equitably
$6$-colorable.

\medskip{}

\noindent 3. One could consider \emph{edge-$k$-critical} graphs with respect to
equitable $k$-coloring, that is, the graphs that are not equitably $k$-colorable
but after deleting any edge become equitably $k$-colorable. One could ask how few
edges may have such $n$-vertex graphs. For $k=3$, they may have asymptotically
$n$ edges, but for $k=4$ possibly this amount 
 is around $5n/4$.

\medskip{}
{ {\bf Acknowledgment.} We thank both referees for their helpful comments.
}

\end{document}